\documentclass[11pt]{article}
\usepackage[utf8]{inputenc}
\usepackage{amsmath}
\usepackage{amssymb}
\usepackage[english]{babel}
\usepackage{hyperref}
\usepackage{apacite}
\usepackage{ragged2e}
\usepackage{enumerate}
\usepackage{authblk}
\usepackage{graphicx}
\usepackage{geometry}
\usepackage{setspace}
\hypersetup{colorlinks=true,linkcolor=blue,anchorcolor=blue,citecolor=blue}
\newtheorem{definition}{Definition}

\newtheorem{theorem}{Theorem}
\newtheorem{lemma}{Lemma}
\newtheorem{proposition}{Proposition}

\newtheorem{assumption}{Assumption}
\newtheorem{corollary}{Corollary}
\title{Bernstein-type Inequalities and Nonparametric Estimation under Near-Epoch Dependence}
\author{Zihao Yuan\thanks{Faculty of Mathematics, Ruhr University Bochum. Financial support from the Deutsche Forschungsgemeinschaft (DFG) TRR 391 Spatio-temporal Statistics for the Transition of Energy and Transport, project number 520388526, is acknowledged.}, Martin Spindler\thanks{Faculty of Business Administration, University of Hamburg. Financial support from the Deutsche Forschungsgemeinschaft (DFG), project number 448795504, is acknowledged.}}
\date{}
\begin{document}

\maketitle
\begin{abstract}
The main contributions of this paper are twofold. First, we derive Bernstein-type inequalities for irregularly spaced data under near-epoch dependent (NED) conditions and deterministic domain-expanding-infill (DEI) asymptotics. By introducing the concept of “effective dimension” to describe the geometric structure of sampled locations, we illustrate – unlike previous research – that the sharpness of these inequalities is affected by this effective dimension. To our knowledge, ours is the first study to report this phenomenon and show Bernstein-type inequalities under deterministic DEI asymptotics. This work represents a direct generalization of the work of \citeA{xu2018sieve}, thus marking an important contribution to the topic. As a corollary, we derive a Bernstein-type inequality for irregularly spaced $\alpha$-mixing random fields under DEI asymptotics. Our second contribution is to apply these inequalities to explore the attainability of optimal convergence rates for the local linear conditional mean estimator under algebraic NED conditions. Our results illustrate how the effective dimension affects assumptions of dependence. This finding refines the results of \citeA{jenish2012nonparametric} and extends the work of \citeA{hansen2008uniform}, \citeA{vogt2012nonparametric}, \citeA{chen2015optimal} and \citeA{li2012local}.
\end{abstract}
\section{Introduction}
\par In time series and spatial analysis, we often encounter different forms of dependence. A natural approach to address this complexity is to introduce a measure of dependence that is broad enough to encompass a large class of time series or spatial models as special cases. The concept of \emph{near-epoch dependence} (NED) was first proposed by \citeA{ibragimov1962some} to model dependence in regularly spaced time series and was later expanded upon by \citeA{billingsley2013convergence} and \citeA{mcleish1974dependent,mcleish1975maximal}. The NED condition is highly general, accommodating sequences that are not identically distributed and encompassing many important processes that do not satisfy mixing conditions. Examples include autoregressive (AR) and infinite moving average (MA($\infty$)) processes under very general conditions, as well as certain dynamic simultaneous equation models (cf. \citeA{bierens2012robust}, Chapter 5 and \citeA{gallant2009nonlinear},  pp 481, 502, 539). In contrast, the $\alpha$-mixing property may not hold under relatively common conditions. For example, \citeA{gorodetskii1978strong} demonstrated that $\alpha$-mixing fails in the case of continuously distributed innovations when the coefficients of a linear process do not decay quickly enough. \citeA{andrews1985nearly} also proved that the AR(1) process with independent Bernoulli innovations does not satisfy $\alpha$-mixing conditions. However, in both cases, the NED condition is satisfied, as established in Proposition 1 of \citeA{jenish2012nonparametric}. Several other regularly spaced time series models satisfy the NED condition. For example, \citeA{hansen1991garch} proved that the GARCH(1,1) process is NED. Furthermore, \citeA{li2012local} showed in Sections 4.3, 4.5 and 4.6 of their paper that the AR(1)-NARCH(1,1), multivariate MA($\infty$) and semiparametric ARCH($\infty$) proposed by \cite{linton2005estimating} also satisfy NED conditions. 
\par With the rise of spatial econometrics, irregularly spaced NED random fields have become increasingly important in many nonlinear spatial processes. \citeA{jenish2012spatial} demonstrated that both linear autoregressive random fields and nonlinear infinite moving average random fields satisfy NED conditions.  \citeA{jenish2012nonparametric} further showed that, under mild moment and Lipschitz conditions, the nonlinear autoregressive random field is also NED, including many widely used spatial auto-regressive (SAR) models as special cases. More recently, \citeA{xu2015spatial} studied the NED properties of SAR models involving a nonlinear transformation of the dependent variable, while \citeA{xu2015maximum} found that, under very mild conditions on the disturbance and weight matrix, the dependent variable in the SAR Tobit model constitutes an $L^2$-NED random field with respect to the disturbance. 

\par More recently, \citeA{xu2018sieve} pointed out that deriving an exponential inequality for irregularly spaced NED random fields is vital for establishing the consistency of sieve-based estimators in many spatial regression models. Even in the case of kernel smoothers, a key step in much of the literature that employs the method of constructing an ``approximate'' process is the application of a Bernstein-type inequality adapted to the dependence conditions of the “approximate” process. Identifying exponential inequalities that directly apply to NED processes, particularly irregularly spaced NED random fields, would greatly simplify such analyses. Lemma 5.1 in \citeA{gerstenberger2018robust} provides a probabilistic inequality for NED processes under a very weak dependence assumption. However, this lemma is limited to data collected from a regularly spaced time series and assumes that the approximate process satisfies $\beta$-mixing conditions, which is insufficiently general to cover our setting. The first major progress on this problem was made by \citeA{xu2018sieve}, who employed the combinatorial techniques of \citeA{doukhan1999new} to establish an exponential inequality under geometric NED conditions of order $e$ (see Definition \ref{Def 1}). Using this result, they proved the consistency of a sieve maximum likelihood estimator for a spatial auto-regressive Tobit model.

\par Another important consideration is the role of infill asymptotics. All of the aforementioned research relies on pure domain-expanding (DE) asymptotics, which ensures that there is always a strictly positive lower bound on the distance between any two locations and this lower bound is independent of sample size. However, as noted by \citeA{lu2014nonparametric}, many important examples inherently require domain-expanding and infill (DEI) asymptotics. Indeed, when selecting locations to observe data, there is often a preference for specific regions, leading to a higher density of sampled locations in those regions. However, this increased density introduces additional complexity due to the irregularity of locations inherent in this approach. Furthermore, the higher density of sampled locations in certain regions imposes stricter dependence conditions, further complicating the theoretical analysis.


\par DEI asymptotics typically rely on one of two assumptions for the sampling of locations: stochastic sampling or deterministic sampling. Under stochastic sampling, the set of selected locations is viewed as a linear transformation of a null set. The key assumption is that for any given sample size $N$, there exists a density function whose support is a finite region (usually a rectangle) such that the null set can be interpreted as an i.i.d. sample generated according to this density. A clear advantage of this approach is that it accommodates non-uniform densities, allowing for selection preferences. However, because this approach treats locations as non-degenerate random vectors, it generally excludes cases in which restrictions on the selection of locations are non-random. For example, spatial locations may be required to be concentrated along multiple lines. Recent research based on this assumption includes \citeA{lahiri2003central}'s central limit theorem for weighted sums and the exploration of the asymptotic properties of frequency-domain empirical likelihood and resampling methods by \citeA{Lahri} and \citeA{lahiri2006resampling}. Additionally, \citeA{kurisu2022nonparametric} and \citeA{kurisu2022local}  have investigated the convergence rates and asymptotic normality of Nadayara-Watson-type estimators for mean and variance functions. The second assumption, deterministic sampling, can easily accommodate non-random restrictions. \citeA{lu2014nonparametric} examined the pointwise asymptotic properties of kernel estimators for marginal and joint density functions under the assumption that infill asymptotics occur uniformly across all locations. Based on this, they introduced an intensity function that serves a similar role to the density function in stochastic sampling. However, this approach does not accommodate cases in which infill asymptotics are restricted to a subset of finite area.

\par Unfortunately, the existing body of research does not extend to the derivation of exponential inequalities under DEI asymptotics and conditions of weak dependence. Additionally, the literature appears to overlook another key aspect, which we refer to in this paper as the effective dimension (see Definition \ref{Definition}). This is puzzling because it challenges a common assumption in spatial data analysis (including all of the previously mentioned research on spatial data) – namely, that increasing $d$ always leads to stronger dependence and more restrictive assumptions provided that the locations lie in $\mathbf{R}^{d}$. However, we present a common counterexample: When $d = 2$ and locations are selected exclusively along the X-axis in $\mathbf{R}^{2}$, the selected locations form a linear-like pattern similar to time series data. Hence, the upper bound of the exponential inequality of averaged partial sums should remain comparable to that in the time series case, suggesting that the higher d does not inherently necessitate stronger dependence conditions or more restrictive assumptions.

\par The contributions of this paper are twofold and can be summarized as follows: First, we derive Bernstein-type inequalities for irregularly spaced data under NED conditions and DEI asymptotics. As a corollary, we also derive a Bernstein-type inequality for irregularly spaced data under $\alpha$-mixing conditions and DEI asymptotics. Importantly, all of our inequalities adapt to the effective dimension (see Definition \ref{Definition}) of the index set, allowing us to obtain sharper upper bounds when the effective dimension is significantly smaller than the dimension of the index set. To our knowledge, ours is the first paper to report this phenomenon. Additionally, apart from the effective dimension, the assumptions regarding DEI asymptotics in this paper are relatively weak and general, requiring only that as the sample size increases infinitely, there exists a sequence of location pairs whose intervening distance vanishes asymptotically. This includes the case in which the set of selected locations comprises a union of two disjoint subsets, one following pure DE asymptotics and the other following DEI asymptotics. Our second contribution is to revisit the attainability of the optimal uniform convergence rate for the local linear estimator of the conditional mean. Unlike nearly all previous studies on this topic, our results highlight how the effective dimension (see Definition \ref{Definition}) affects the dependence assumptions required to achieve optimal rates.


\par The remainder of this paper is structured as follows: Section 2.1 introduces the assumptions related to selecting locations, including DEI asymptotics. Section 2.2 defines NED conditions and explores some of their useful properties. Sections 3 and 4 are dedicated to an in-depth discussion of the two main contributions outlined above. Lastly, the proofs are presented in the Appendix.



\section{Irregularly Spaced NED Random Fields}
This section introduces our definition of irregularly spaced NED random fields. In Section 2.1, we illustrate our assumptions governing the selection of locations. Section 2.2 defines NED conditions and explores some of their interesting properties.
\subsection{Selection of Locations}
\begin{assumption}
\label{Asp 1}
Individual units are located at $\Gamma_{N}\subset\mathbf{R}^{d}$, where $d$ is a positive integer. Here $\Gamma_{N}:=\{s_{iN}: i=1,2,...,N\}$ is the set of selected locations, which we assume are selected in a deterministic manner.
\end{assumption}
\par Assumption \ref{Asp 1} excludes the case of stochastic sampling. As discussed in Section 1, we focus exclusively on deterministic sampling in this paper because it readily incorporates some subjective preferences or other restrictions on location selection.
\begin{assumption}
\label{Asp 2}
For any two distinct selected locations $s_{iN},s_{jN}\in \Gamma_{N}$, we assume that $||s_{iN}-s_{jN}||>d_{0}$ where $||\cdot||$ denotes the Euclidean norm and $d_{0}>0$ for each $N$. Furthermore, we assume that the following two conditions about $d_0$ hold simultaneously: (1) $\lim_{N}Nd_{0}^{d}\nearrow+\infty$; (2) $d_0$ is non-increasing as sample size $N\nearrow +\infty$.   
\end{assumption}
The condition $\lim_{N}Nd_{0}^{d}\nearrow+\infty$ ensures that, as $N\nearrow +\infty$, the set of selected locations can never be contained within a bounded subset of $\mathbf{R}^d$, necessitating an asymptotic expansion of the domain. In this paper, we say that the index set $\Gamma_N$ (i.e., the set of selected locations) follows \textbf{domain-expanding (DE) asymptotics} if $\lim_{N}Nd_{0}^{d}\nearrow+\infty$ and $\inf_{N}d_0>0$ hold. If, instead, $\lim_{N}Nd_{0}^{d}\nearrow+\infty$ and $\lim_Nd_0\searrow 0$ hold, we say $\Gamma_N$ has \textbf{domain-expanding and infill (DEI) asymptotics}. 

\par Provided that $i \neq j$ indicates $s_{iN}\neq s_{jN}$, an obvious equation for $d_0$ is: 
\begin{align*}
    d_0=\min_{1\leq i\leq N}\min\{||s_{jN}-s_{iN}||: 1\leq j\neq i\leq N\}=:\min_{1\leq i\leq N}\delta_{i}.
\end{align*} 
Hence, the “infill” part of our DEI asymptotics requires
\begin{align*}
    \lim_{N}\min_{1\leq i\leq N}\delta_{i}=0.
\end{align*}
To our knowledge, \citeA{lu2014nonparametric} were the first to provide a clear and practical definition of deterministic DEI asymptotics. Based on formula (2.3) in their paper, but using our notation above, the infill part of their DEI asymptotics requires
\begin{align*}
    \lim_{N}\max_{1\leq i\leq N}\delta_{i}=0.
\end{align*}
Thus, our formulation of infill asymptotics is more general than that of \citeA{lu2014nonparametric}, whose definition can be considered a uniform counterpart to our assumption. Unlike their approach, which assumes dense observation across all locations, our concept of DEI asymptotics allows for selective denseness in specific regions of interest. This flexibility makes our approach more adaptable to various practical scenarios and more accurately mirrors real-world data collection processes. 
\par Regardless of whether infill asymptotics are incorporated, the following proposition plays a key role in our proof of Theorem \ref{Theorem 1} and its corollaries. It enables us to approximate complex, irregularly spaced settings using regularly spaced and “group variable” approaches. Moreover, this proposition represents a modification of Lemma A.1 from \citeA{jenish2009central}\footnote{We note an ambiguity in the calculation of the proof of Lemma A.1 in \citeA{jenish2009central}. To address this, we provide an alternative proof of the proposition. The only difference lies in the constant factors; however the two results are asymptotically equivalent.}
\begin{proposition}
\label{prop 1}
Suppose Assumptions \ref{Asp 1} and \ref{Asp 2} hold. Then,
\begin{enumerate}
    \item[\rm{(\textbf{i})}] by denoting $I(s,h)$ as a closed cube centered at location $s\in \Gamma$ with the length of each edge as $h$, we have $|I(s,\frac{\sqrt{2}d_0}{2})\cap \Gamma_{N}|\leq 1$, where $|\cdot|$ denotes the cardinality of a set.
    \item[\rm{(\textbf{ii})}] there exists a constant $C<+\infty$ such that for $h\geq 1$
    $$\sup_{s\in \Gamma}|I(s,h)\cap \Gamma_{N}|\leq Ch^{d},\ C=(2\sqrt{2}/d_0)^d.$$
\end{enumerate}
\end{proposition}
Point (i) of Proposition \ref{prop 1} reflects, how “small” a cube must be, given the minimum distance $d_0$, to ensure that it contains at most one location. Point (ii) gives an upper bound of the number of locations contained in a given cube. Before presenting Assumption \ref{Asp 3}, we introduce the following definition of the effective dimension of the index set  $\Gamma_{N}$. 
\begin{definition}(\textbf{Effective Dimension})
    \label{Definition} 
    {For a positive integer $l$, let $R_{l}(r_{l})$ denote a rectangle in $\mathbf{R}^{d}$ such that $r_{l}$ out of $d$ edges are infinitely long while the length of the remaining $d-r_{l}$ edges is fixed. Based on the set of selected locations $\Gamma_{N}$ defined in Assumption \ref{Asp 1}, and denoting $\mathbf{N}$ as the set of natural numbers, we define $\Gamma(d')$ as a “\textbf{$d'$-dimensional candidate set}” if there exists a universal finite positive integer $M$ and an $M$-tuple  $\mathbf{r}:=(r_{1},r_{2},\dots, r_{M})\in\mathbf{N}^{M}$ such that
    \begin{itemize}
        \item [1.] $\max_{1\leq l\leq M}r_l= d'$,
        \item [2.] $\Gamma_{N}\subset \bigcup_{l=1}^{M} R_{l}(r_{l})=:\Gamma(d')$ holds for any $N$.
    \end{itemize}
    We then say that set $\Gamma_{N}$ has an \textbf{effective dimension} $d_2$ if
    \begin{align*}
        d_2=\min\{ 1\leq d'\leq d:\ \exists\  \Gamma(d')\ \text{serving as a $d'$-dimensional candidate set}\}.
    \end{align*}
    }
\end{definition}
\par
\begin{figure}
    \begin{center}
        \includegraphics[width=0.6\linewidth]{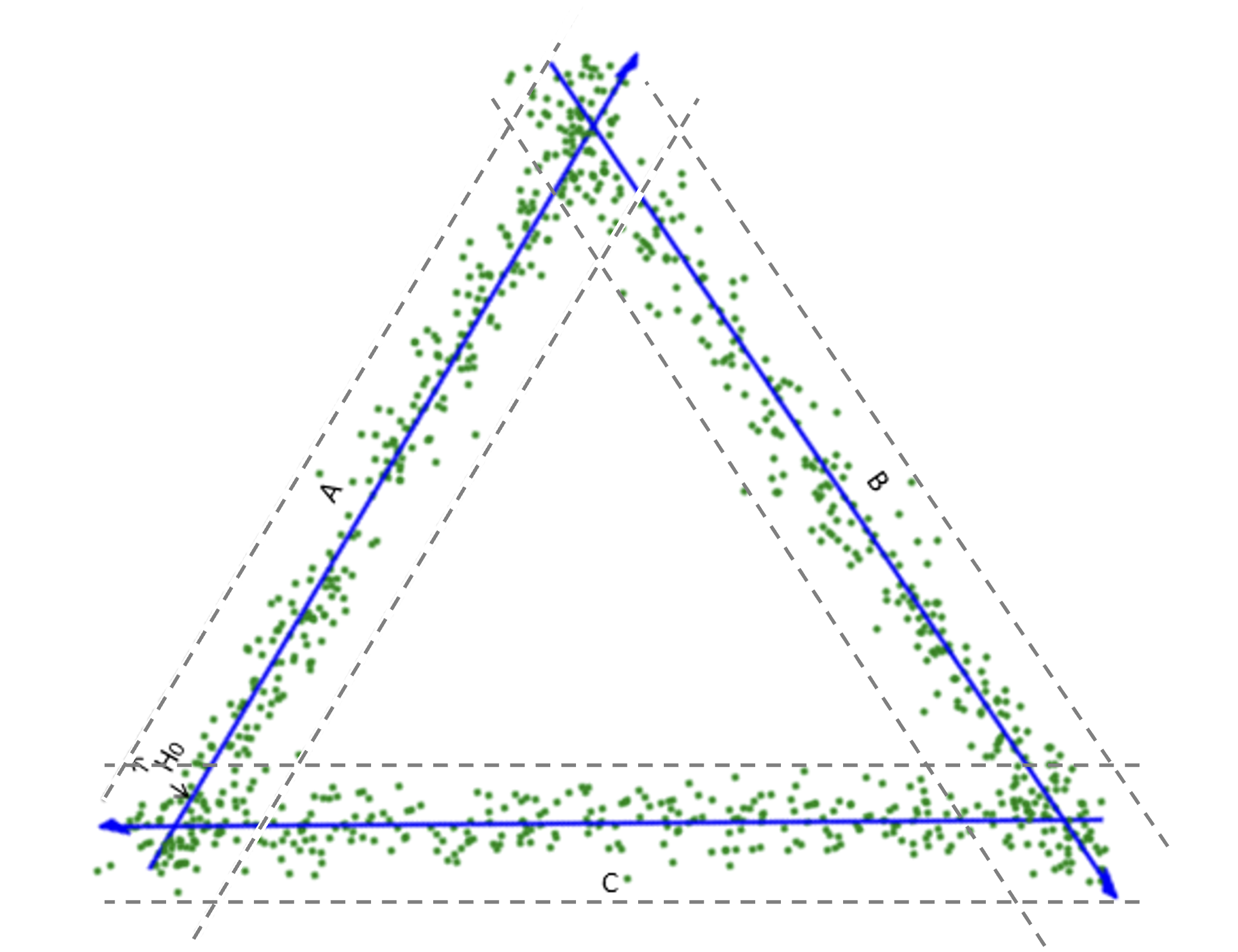}
    \end{center}
    \caption{\textit{Green dots denote the selected locations. These locations are selected from the candidate set, which is the union of the $H_{0}$ blow-ups of lines A, B, and C, where $H_{0}$ is universally finite.}}
    \label{Figure 1}
\end{figure}

 Definition \ref{Definition} states that when $\Gamma_N$ has an effective dimension $d_2$ strictly smaller than $d$, the selected locations can be covered by a collection of finite blow-ups\footnote{Given a set $A$ in a metric space $(E,d)$, for some $h>0$, we define its $h$-blow-up as: $A^{+h}=\{x\in E: d(x,A)\leq h\}$. Naturally, $A^{+h}$ is considered a finite blow-up if $h<+\infty$.} of $d_2$-dimensional Euclidean subspaces. The definition allows for overlapping rectangles in this construction. If we ignore the possibility of a lower effective dimension by setting $d_2=d$, we can directly use $\mathbf{R}^{d}$ as the rectangle described in Definition \ref{Definition}. In this case, $M$ is equal to $1$, which is trivial. Figure \ref{Figure 1} provides a non-trivial example of the existence of a lower effective dimension, i.e. $d_2<d$. 
 In Figure \ref{Figure 1}, the candidate set is represented as the union of three $\mathbf{H}_{0}$ blow-ups of three infinitely long lines in $\mathbf{R}^{2}$, resulting in $M=3$ while $d_2=1$. This implies that the dependence among spatial data, characterized by the distribution of locations in Figure \ref{Figure 1}, is asymptotically no stronger than that of three groups of ``time-series-like" spatial data, because for each group, only one direction is relevant.
 \par Moreover, it is important to recall that, based on the fundamental properties of the sum operation for real-valued random variables, we can always express the total partial sum of the  collected spatial data as the accumulation of multiple sub-group partial sums. By defining these sub-groups as the three blow-ups in Figure 1, we decompose a single “sum of spatial data” into three “partial sums of time-series-like spatial data”. Because the central task of this paper is to derive exponential inequalities, we draw upon sub-additivity of measures, to obtain:
 \begin{align*}
     \mathbb{P}\Big(\Big|\sum_{s_{iN}\in\Gamma_N}X_{s_{iN}}\Big|>t\Big)\leq \sum_{k=A,B,C}\mathbb{P}\Big(\Big|\sum_{s_{iN}\in \mathbf{B}_{k}}X_{s_{iN}}\Big|>t\Big),
 \end{align*}
 where $\mathbf{B}_{k}$ denotes the $H_0$ blow-up of line $k$ for  $k=A,B,C$. This also explains why we restrict $M$ to be universally finite in Definition \ref{Definition}. Indeed, without additional preconditions, taking advantage of the lower effective dimension of the index set heavily relies on the sub-additivity property mentioned above. If $M$ were to grow with the sample size, the sharpness of the  inequalities would be compromised. Thus, we explicitly exclude this situation in our definition. 
\par According to Definition \ref{Definition}, for any set of observed spatial data whose locations satisfy Assumptions 1 and 2, a potentially lower effective dimension $d_2$ is always well defined and exists.  Furthermore, due to the finiteness of $M$, the sub-additivity of measures, and Theorem \ref{Theorem *}, we know that the sharpness of inequalities is determined only by $d_2$ rather than $d$. A practical and efficient way to identify $d_2$ is through prior knowledge of the sample collection. For example, if the spatial data set is known to have been collected within some finite blow-up of several given rivers or mountain paths, similar to our previous discussion of Figure 1, we can immediately conclude that $d_2=1$ while $d=2$. It is also of interest to develop a non-asymptotic definition of effective dimension that is fully consistent with Definition \ref{Definition}. 
\par \par
\begin{figure}
    \begin{center}
        \includegraphics[width=0.6\linewidth]{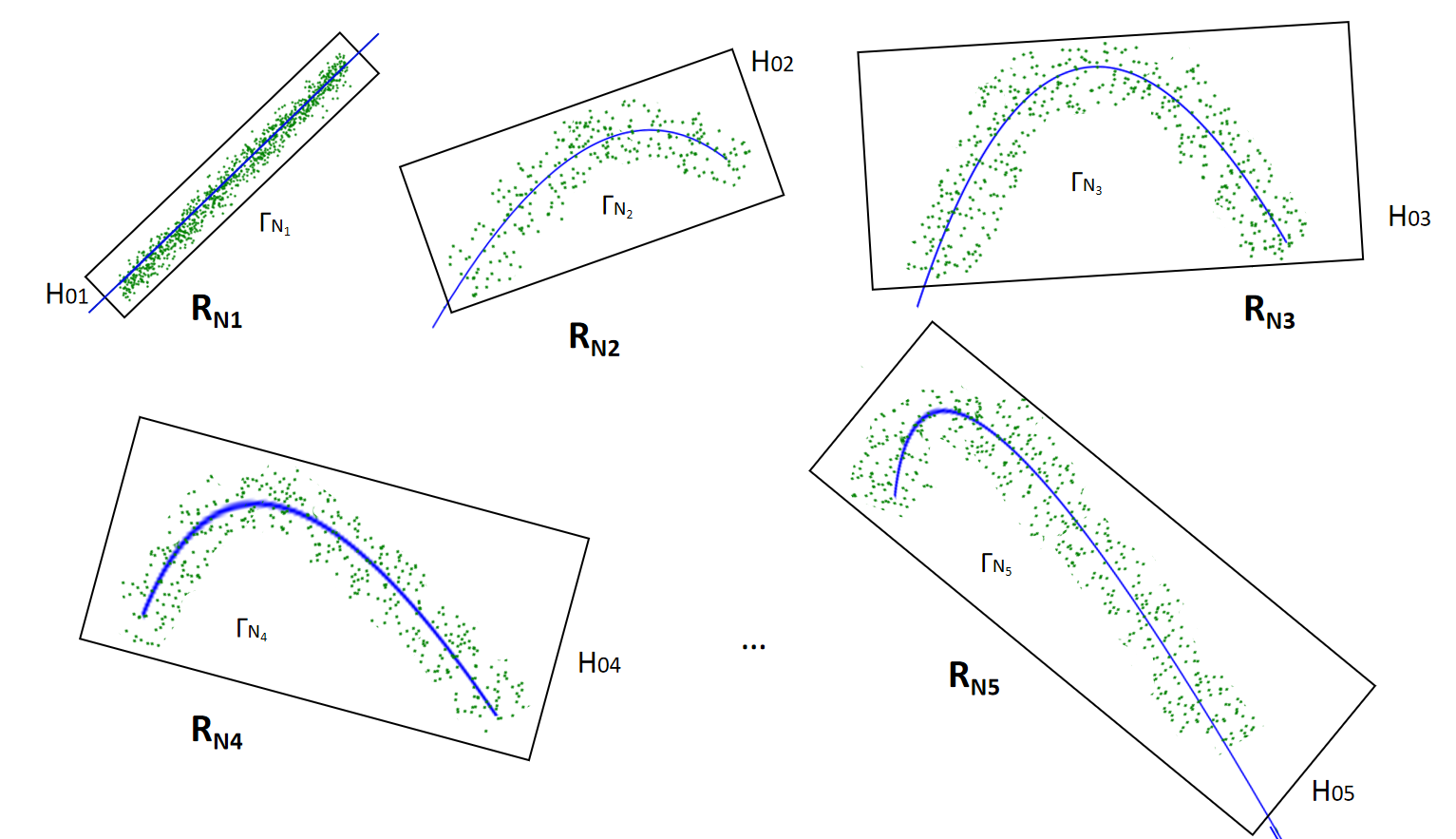}
    \end{center}
    \caption{\it The $N_k$s denote the sample size and satisfy $N_k\leq N_{k+1}$. Based on the definition of $R_{N_k}$ and $\Gamma_{N_k}$, we have $\Gamma_{N_k}\subset \Gamma_{N_{k+1}}$. Thus, Figure 2 illustrates the variation of  $R_{N_k}$ as the sample size grows. Assume that the selected locations (green dots) lie only within a finite blow-up of the blue line. If the change in the blue line in some given direction is uniformly bounded, we can regard $\max_{k\in\mathbf{N}}H_{0k}$ as a finite positive number, since an infinitely long band with finite bandwidth is large enough to contain all locations (green dots) in $\Gamma_N$. Meanwhile, the change in the blue line cannot be a graph, such as the Archimedean spiral, which is therefore excluded.}
    \label{Figure 2}
\end{figure}

 Because Definition \ref{Definition} strongly relies on the fact that $\Gamma_N$ is contained by some candidate set $\Gamma(d')$, we now introduce the following necessary and sufficient condition for $\Gamma_N\subset \Gamma(d')$:
\begin{proposition}
\label{prop +}
Building upon the concept of $\Gamma(d')$ introduced in Definition \ref{Definition}, a necessary and sufficient condition for $\Gamma_{N}\subset \Gamma(d')$ to hold uniformly over $N$ is as follows: Given $\Gamma_{N}$, there exists a collection of rectangles, denoted as $\{R_{lN}\subset \mathbf{R}^{d}: 1\leq l\leq M_{N},\ N\in\mathbf{N}^{+}\}$, where ${M_{N}}\in \mathbf{N}^{+}$, satisfying the following conditions: 
\begin{itemize}
    \item[1] $\Gamma_{N}\subset \bigcup_{l=1}^{M_{N}}R_{lN}$ holds for each $N$ and $\max_{N}M_{N}\leq M$.
    \item[2] For every $N$ and $l$, the volume of $R_{lN}$ is finite but tends to diverge as $N$ approaches $+\infty$.
    \item[3] For each $1\leq l\leq M_{N}$, there exists a $H_{l}>0$ independent of $N$ such that only $d_{l}$ out of $d$ edges could surpass a length of $H_{l}$ and approaches $+\infty$. Furthermore, $\max_{1\leq l\leq M_{N}}H_{l}<+\infty$ and $\max_{1\leq l\leq M_{N}}d_{l}= d'$, for some $1\leq d'\leq d$.  
\end{itemize}
\end{proposition} 
     Figure \ref{Figure 2} gives an example with $M=1$ and $d=2$, illustrating how the volume of rectangle $R_{N}$ becomes larger as $N$ increases. This is due to the assumption that $Nd_{0}^{d}\nearrow +\infty$. However, if the selected locations (green dots) show path dependence -- meaning that the maximum distance from the curve to each location (blue curves) remains uniformly bounded –- a finite blow-up of the blue curve can encompass all selected locations. Thus, for any given $N$, a uniform upper bound exists for the length of the shorter edges in each $R_{N}$, satisfying conditions 2 and 3 in Proposition \ref{prop +}.  
    \par The array $\{R_{lN}\}$ introduced in Proposition \ref{prop +} then yields the following non-asymptotic definition of effective dimension, which is equivalent to Definition \ref{Definition}.
\begin{definition}
    \label{Definition +} (Non-Asymptotic \textbf{Effective Dimension})
    For a given $\Gamma_{N}$ satisfying Assumptions \ref{Asp 1} and \ref{Asp 2}, we say it has an effective dimension $d_2$ if 
    \begin{align*}
        d_{2}=\min_{\{R_{lN}\}}\max_{1\leq l\leq M_{N}}d_{l},
    \end{align*}
    where the $\min$ is taken over all $\{R_{lN}\}$ that satisfy conditions 1, 2, and 3 in Proposition \ref{prop +}.
\end{definition}
Then, according to the necessity and sufficiency condition established in Proposition \ref{prop +}, we can conclude that Definition \ref{Definition +} is equal to Definition \ref{Definition}.
\begin{assumption}
\label{Asp 3}
 According to Definitions \ref{Definition} and \ref{Definition +}, we assume that $\Gamma_{N}$ has an effective dimension $d_2$, where $1\leq d_2\leq d$ and that $M=1$.      
\end{assumption} 
\par
As per Definition \ref{Definition}, Assumption \ref{Asp 3} suggests that the candidate set can be regarded as a finite blow-up of $\mathbf{R}^{d_2}$ within $\mathbf{R}^{d}$. For simplicity, we assume that the length of this blow-up is $H_{0}$\footnote{This implies that the “candidate set” resembles a rectangle, with $d_2$ out of $d$ edges being infinitely long, while the remaining edges have a length of $H_{0}$, where $H_{0}$ is a universal positive constant.}. Thus, based on the equivalence between Definitions \ref{Definition} and \ref{Definition +}, for any given sample size $N$, a finite rectangle-like subset of this blow-up, denoted as $R_{N}$, suffices to encompass all selected locations. For simplicity, we set $R_N$ as a truncation of $\mathbf{R}^{d_2}$, meaning that $d_2$ out of $d$ edges of $R_N$ diverge to $+\infty$ as $N\nearrow +\infty$, while the remaining edges have a length of $H_{0}$.
\begin{assumption}
\label{Asp 4}
For any given sample size $N$ and the rectangle $R_N$ mentioned above, we assume $$Vol(R_N)\leq KB_{N},$$ for some $K>0$ and $B_{N}>0$ holding for each $N$ with $\lim (B_{N}/Nd_{0}^{d})>0$. In this paper, we primarily consider cases where $B_{N}=Nd_{0}^{d}$ or $B_N=N$, but our theorems hold valid for every given $B_N>0$. 
\end{assumption}
\par Generally speaking, Assumption 4 ensures that the distribution of locations is not excessively sparse. More importantly, it plays two key roles: First, it establishes a link between sample size and the number of blocks constructed in our proof of inequalities. Second, it facilitates the application of our “group-variable” approach to handle the irregularity of locations.

\par To clarify this strategy, recall that we can partition $R_N$ into numerous, evenly spaced small blocks. Proposition \ref{prop 1} enables us to establish an upper bound on the number of selected locations within each block. By treating the sum of random variables within each block as a new random variable (henceforth referred to as “group variables”), we express the total sum of the original random variables as the sum of these “group variables”, which are regularly spaced in $R_N$. This effectively transforms an irregularly spaced setting into a regularly spaced one. This transformation allows Bernstein's well-established blocking technique, which was originally developed for regularly spaced data, to be adapted to the more complex case of  irregularly spaced data. Importantly, this construction holds for any given $d_{0}>0$, thus accommodating infill asymptotics in a very general manner. However, in the context of blocking techniques, the number of grouped variables  serves as the effective sample size. Hence, a key challenge is to establish the connection between the number of “grouped variables” and the true sample size $N$. This relationship is formalized through Assumption \ref{Asp 4}.

\begin{assumption}
    \label{Asp 4+}
 We assume $\min_{1\leq j\leq d_{2}}H_{k_{j}}^{N}\geq \beta(N)$, where $\beta(N)>0$ for any $N$ and $\lim_{N}\beta(N)\nearrow +\infty$.
\end{assumption}
 \par Generally speaking, if we understand Assumption \ref{Asp 4} as a restriction on the upper bound of the growth condition of maximum distance between any two selected locations, Assumption \ref{Asp 4+} serves as a lower bound. In the context of regularly spaced data, such an assumption is common (e.g., \citeA{liebscher1996strong} and \citeA{yao2003exponential}). Although our proof still relies on Bernstein's blocking technique, the irregularity of the selected locations, unlike in regularly spaced data, forces us to apply blocking to the sampling region $R_{N}$. This introduces challenges in linking the block size to the sample size. Fortunately, Assumptions \ref{Asp 4+} and \ref{Asp 4} allow us to establish this connection in a straightforward manner. Additionally, in Sections 3 and 4, we assume that the growth condition of $\beta(N)$ is at least $Cp_N$. This assumption is mainly made for  technical convenience, preventing the blocks we construct in our proof from overestimating the number of locations in each block. However, even if $\beta(N)$ grows more slowly than the block size, our proof strategy remains valid -- the only difference being --  that the calculations become considerably more complex. This occurs because each block is no longer a standard cube but instead takes the shape of a rectangle, requiring the introduction of two parameters, $p_N$ and $p'_N$, to measure the block size, where $p'_N=\mathcal{O}(\beta(N))$. For simplicity, in this paper, we only consider cases where $\beta(N)$ grows at least as fast as the block size.
\subsection{NED Conditions}
\par For any random vector $v\in \mathbf{R}^{D}$, $||v||_{p}=[E||v||^{p}]^{\frac{1}{p}}$ denotes the $L^{p}$-norm of the Euclidean norm of vector $v$ with respect to a probabilistic measure. For any subset $A\subset \mathbf{R}^{d}$, we denote its cardinality as $|A|$. Based on Assumption \ref{Asp 1}, we introduce the definitions and properties of $\mathbf{R}^{D}$-valued ($D\geq 1$) NED random fields as follows. 
\begin{definition}
\label{Def 1}  
Let $Z=\{Z_{s_{iN}}:s_{iN}\in \Gamma_{N}\}$ be an $\mathbf{R}^D$-valued array and let $\epsilon=\{\epsilon_{s_{iN}}: s_{iN}\in \Gamma_{N}\}$ be another array, where $\Gamma_{N}$ satisfies Assumption \ref{Asp 1}. Define $\mathcal{F}_{i,N}(r)$ as the $\sigma$-field generated by the random vectors $\epsilon_{s_{jN}}'s$ satisfying $||s_{jN}-s_{iN}||\leq r$. Then, the array $Z$ is said to be $L^p$-near-epoch dependent ($L^p$-NED) on $\epsilon$ if, for some $p>0$, we have
\begin{equation}
    \label{eq:1}
    ||Z_{s_{iN}}-E[Z_{s_{iN}}|\mathcal{F}_{iN}(r)]||_{p}\leq \alpha_{s_{iN}}\psi_{Z}(r),
\end{equation}
where $\{\alpha_{s_{iN}}\}$ is an array of finite positive constants associated with locations $\{s_{iN}\}$ and function $\psi_{Z}(r)$ is a positive-valued non-increasing function that satisfies $\lim_{r\rightarrow +\infty}\psi_{Z}(r)=0$. We say $Z$ is \textbf{geometrically } $L^p$-NED on $\epsilon$ of order $(\nu,b,\gamma)$, if $\psi_{Z}(r)\leq \nu^{-br^{\gamma}}$, for some $\nu>1,\ b,\ \gamma>0$. Similarly, Z is algebraically $L^p$-NED on $\epsilon$ of order $\gamma$ if
$$ \psi_{Z}(r)\leq r^{-\gamma}, \text{for some}\ \gamma>0.$$
\end{definition}
Definition \ref{Def 1} is borrowed directly from \citeA{jenish2012spatial}. As we noted in Section 1, a widely used setting -- also adopted in this paper -- is to assume that the random field $\epsilon$ satisfied a strong mixing conditions. To lay the groundwork for subsequent sections, we begin by revisiting the definition of a strong mixing random field.
\begin{definition}
\label{Def 2}
Given that $\Gamma_N$ satisfies Assumption \ref{Asp 1}, we consider that an array $\epsilon:=\{\epsilon_{s_{iN}}: s_{iN}\in\Gamma_{N}\}$ satisfies $\alpha$-mixing conditions if, for any subsets $ U,\ V\subset \Gamma_{N}$ such that $\max\{|U|,|V|\}<+\infty$, the following inequality holds:
\begin{align*}
    &\alpha(\sigma(U),\sigma(V))\\
    &:=\sup\{|\mathbb{P}(A \cap B)-\mathbb{P}(A)\mathbb{P}(B)|, A\in \sigma(U), B\in \sigma(V)\}\\
    &\leq \phi(|U|,|V|)\psi_{\epsilon}(\rho(U,V)),
\end{align*}
where $\sigma(U)$ denotes the sigma-algebra generated by $\epsilon_{s_{iN}}$ whose coordinates belong to the set $U$. Function $\phi:\mathbf{N}^2\rightarrow \mathbf{R}^{+}$ is symmetric and increasing in both of its arguments; $\psi_{\epsilon}(\cdot)$ has the same properties as the function $\psi_{Z}(\cdot)$ defined in Definition \ref{Def 1}; and $\rho(U,V)=\inf\{||u-v||,u\in U, v\in V\}$, where $||\cdot||$ denotes the Euclidean norm.
\end{definition}
Based on Definitions 3 and 4, we can now obtain the following propositions which will later be useful in deriving concentration inequalities for arrays of random variables under NED conditions.
\begin{proposition}
\label{prop 2}
Given an array $Z=\{Z_{s_{iN}}: s_{iN}\in \Gamma_{N}\}$ and $\epsilon=\{\epsilon_{s_{iN}}: s_{iN}\in\Gamma_{N}\}$, where $\Gamma_N$ satisfies Assumption 1 and $Z_{s_{iN}}=(Z^{1}_{s_{iN}},...,Z^{D}_{s_{iN}})\in\mathbf{R}^{D}$, $Z$ is $L^{p}$-NED on $\epsilon$ if and only if for every $k=1,...,D$, $\{Z^{k}_{s_{iN}}:s_{iN}\in\Gamma_{N}\}$ is also $L^p$-NED on $\epsilon$.   
\end{proposition}
\begin{proposition}
\label{prop 3}
Suppose the array $Z$, as  defined in Definition \ref{Def 1}, is $L^{p}$-NED on the random field $\epsilon$. Additionally, for any $s_{iN}\in\Gamma_N$, we have $Z_{s_{iN}}\in\mathcal{Z}\subset \mathbf{R}^D$. Let $f$ be any real-valued measurable function whose domain is $\mathcal{Z}$. If $f$ is coordinate-wise Lipschitz continuous, with $Lip_{k}(f)$ denoting the Lipschitz constant for the $k$-th coordinate, then the transformed array $f(Z):=\{f(Z_{s_{iN}}), s_{iN}\in \Gamma_N\}$ is also $L^p$-NED on $\epsilon$. Furthermore, we have the following bound:
\begin{equation}
\label{eq:2}
||f(Z_{s_{iN}})-E[f(Z_{s_{iN}})|\mathcal{F}_{iN}(r)]||_{p}\leq 2\sum_{k=1}^{D}Lip_{k}(f)\alpha_{s_{iN}}\psi_{Z}(r).
\end{equation}
\end{proposition}
\begin{proposition}
\label{prop 4}
Suppose that the array $Z$ is $L^{p}$-NED on a strongly mixing random field $\epsilon$, which satisfies Definition \ref{Def 2}. Additionally, assume that for any $s_{iN}\in\Gamma_N$, we have  $Z_{s_{iN}}\in\mathcal{Z}\subset \mathbf{R}$. For $\forall\ U,\ V\subset \Gamma_{N}$, such that $\max\{|U|,|V|\}<+\infty$, denote $U=\{s_{iN},i=1,...,|U|\}$, $V=\{t_{jN},j=1,...,|V|\}$. Let $g:\mathcal{Z}^{|U|}\rightarrow R$, $h:\mathcal{Z}^{|V|}\rightarrow R$ be two coordinate wise Lipschitz functions with $Lip(g)$ and $Lip(h)$ denoting their respective Lipschitz constants. If $p\geq 1$, for some $\delta>2-p$ and $r<\rho(U,V)/2$, we have
\begin{align}
&\left |Cov(g(Z_{s_{1N}},...,Z_{s_{|U|N}}),h(Z_{t_{1N}},...,Z_{t_{|V|N}}))\right |\notag\\
&\leq 4Lip(f)Lip(g)\begin{matrix}(\sum_{i=1}^{|U|}\alpha_{s_{iN}})((\sum_{j=1}^{|V|}\alpha_{t_{jN}}))\end{matrix}\psi_{Z}^{2}(r)\notag\\
&+4Lip(g)\begin{matrix}(\sum_{i=1}^{|U|}\end{matrix}\alpha_{s_{iN}})\psi_{Z}(r)||h||_{p}\notag\\
&+4Lip(f)\begin{matrix}(\sum_{j=1}^{|V|}\end{matrix}\alpha_{t_{jN}})\psi_{Z}(r)||g||_{p}\notag\\
&+\phi(|U|,|V|)\psi_{\epsilon}^{\frac{p+\delta-2}{p+\delta}}(\rho(U,V)-2r)||g||_{p+\delta}||h||_{p+\delta}.\label{eq:3}
\end{align}
\end{proposition}
\par Proposition \ref{prop 2} is straightforward and easy to prove using only the basic property of conditional expectations and the relationship between $L^1$ and Euclidean norm in $\mathbf{R}^D$. Proposition \ref{prop 3} serves as the foundation for Proposition \ref{prop 4} and allows us to concentrate solely on real-valued NED processes when the statistic can be described as an empirical process whose index set is Lipschitz class. Proposition \ref{prop 4} serves as a refined version of Lemma A.2 from \citeA{xu2018sieve}. The only difference is that our Equation \eqref{eq:3} relies on the $L^{p}$-norm instead of the $L^{\infty}$-norm. Proposition \ref{prop 2} allows us to apply the concentration inequalities displayed in Section 3 directly to many interesting estimators of regression functions. After all, under mild conditions, most of these can be regarded as Lipschitz functions of the data, such as kernel or kNN smoothers, Hermite polynomials or various other orthogonal series estimators. 
\par {It is commonly understood that the approach described in Proposition \ref{prop 4} plays an important role in establishing upper bounds for sums of random variables that satisfy NED conditions and specific location arrangements.
Theorem 1 in \citeA{hansen2008uniform} addresses a similar  topic in the context of  strong mixing conditions and regularly spaced locations. Thus, the following proposition can be regarded as a result parallel to \citeauthor{hansen2008uniform}'s Theorem 1.
\begin{proposition}
    \label{prop *}
    Suppose $\Lambda$ is a cube in $R^{d}$ defined as follows,
    \begin{align*}
        \Lambda=\prod_{k=1}^{d_2}(0,h_{k}]\times (0,h_0)^{d_1},
    \end{align*}
    where $d_2+d_1=d$ and $\min_{1\leq k\leq d_2}h_{k}\geq h_{0}$. Let $Z=\{Z_{s_{iN}}: s_{iN}\in \Gamma_{N}\}$ be $L^{p}$-NED ($p\geq 1$) on $\epsilon=\{\epsilon_{s_{iN}}: s_{iN}\in \Gamma_{N}\}$, where $Z_{s_{iN}}\in\mathbf{R}$ and $\Gamma_{N}$ is the set of selected locations submitted to Assumptions \ref{Asp 1} and \ref{Asp 2}. For some $\delta>2-p$, assume $\max_{in}||Z_{s_{iN}}||_{2}\leq \sigma$, $\max_{in}||Z_{s_{iN}}||_{2+\delta}\leq \sigma_{2+\delta}$ and $\max_{i\neq j}E|Z_{s_{iN}}Z_{s_{jN}}|+E|Z_{s_{iN}}|E|Z_{s_{jN}}|\leq \Sigma$. Additionally, we assume $\max_{i,N}\alpha_{s_{iN}}\leq A(N)$, where $A(N)$ is positive for each $N$. Then, for any given non-decreasing positive integer sequence $\{\tau_{N}\}$, we have:
    \begin{align*}
     & Var(\sum_{s_{iN}\in\Lambda}Z_{s_{iN}})\leq N(\Lambda)\mathbf{V},\ \text{where}\\
     &\mathbf{V}=\sigma^{2}+2^{d_2}n_{0}^{d_1}\tau_{N}^{d_2}\Sigma+ (2\pi)^{d_2}h_{0}^{d_1}d_{0}^{-d} (\mathbf{T}_{1}+\mathbf{T}_{2}+\mathbf{T}_{3}),\\
     &\mathbf{T}_1= 4A(N)\int_{\frac{\tau_{N}d_{0}}{\sqrt{2}}}^{+\infty}\rho^{d_2-1}\psi_{Z}^{2}(\frac{\rho}{3})d\rho,\ \mathbf{T}_{2} = 8A(N)\sigma\int_{\frac{\tau_{N}d_{0}}{\sqrt{2}}}^{+\infty}\rho^{d_2-1}\psi_{Z}(\frac{\rho}{3})d\rho,\\
     &\mathbf{T}_3 =2^{\tau}\sigma_{2+\delta}^{2}\int_{\frac{\tau_{N}d_{0}}{\sqrt{2}}}^{+\infty}\rho^{d_2-1}\psi_{\epsilon}^{\frac{\delta}{2+\delta}}(\frac{\rho}{3})d\rho, N(\Lambda)= |\Gamma_{N}\cap \Lambda|,
    \end{align*}
    and “$n_0$” is some universal number such that $n_0^d\leq ([h_0/d_0]+1)^d$. More specifically, when $\psi_{Z}(r)=\psi_{\epsilon}(r)=\nu^{-br^{\gamma}}$, $\nu>1$ and $b,\gamma>0$, there exists some universal constant $C^*>0$ independent of $N$ such that
    \begin{align*}
        \sum_{l=1}^{3}\mathbf{T}_l\leq C^{*}(\frac{\tau_{N}d_0}{2\sqrt{2}})^{d_{2}-\gamma}&\left(A(N)\nu^{-2b(\frac{\tau_{N}d_0}{2\sqrt{2}})^{\gamma}}+A(N)\sigma \nu^{-b(\frac{\tau_{N}d_0}{2\sqrt{2}})^{\gamma}} \right.\\
        &\left.+\sigma_{2+\delta}^{2}\nu^{-\frac{b\delta}{(2+\delta)}(\frac{\tau_{N}d_0}{2\sqrt{2}})^{\gamma}}\right):=\mathbf{B1}.
    \end{align*}
    Hence, we denote $\mathbf{V1}=\sigma^{2}+2^{d_2}n_{0}^{d_1}\tau_{N}^{d_{2}}\Sigma+(2\pi)^{d_2}h_{0}^{d_1}d_{0}^{-d}\mathbf{B1}.$ When $\psi_{Z}(r)=\psi_{\epsilon}(r)=r^{-\gamma}$ and $\gamma>\frac{d_{2}(2+\delta)}{\delta}$, there exists some $C^{**}>0$ independent of $N$ such that  
   \begin{align*}
    \sum_{l=1}^{3}\mathbf{T}_l \leq C^{**}&\left(A(N)(\tau_{N}d_{0})^{-(2\gamma-d_{2})}+\sigma A(N)(\tau_{N}d_{0})^{-(\gamma-d_{2})}\right. \\
    &\left.+\sigma_{2+\delta}^{2}(\tau_{N}d_{0})^{-(\frac{\delta\gamma}{2+\delta}-d_{2})}\right):=\mathbf{B2}.
\end{align*}
Similarly, we denote 
$\mathbf{V2}=\sigma^{2}+2^{d_2}n_{0}^{d_1}\tau_{N}^{d_{2}}\Sigma+(2\pi)^{d_2}h_{0}^{d_1}d_{0}^{-d}\mathbf{B2}.$
\end{proposition}
The reason for introducing number sequence $\tau_N$ in Proposition \ref{prop *} is purely technical and largely motivated by the proof of Theorem 1 in \citeA{hansen2008uniform}. The  specific advantages of incorporating $\tau_N$ become clearer in the proof of Theorem 2. It is evident that when both NED and $\alpha$-mixing exhibit algebraic decay, the inclusion of infill asymptotics, particularly as $d_0$ decreases to zero ($d_0\searrow 0$), continues to  influence the sharpness of the variance inequality. The extent of this influence is determined by assumptions on dependence and the interplay between $d_{0}$ and tuning parameters, such as the bandwidth of kernel smoothers. This relationship will be further elucidated in Section 4. However, no matter how carefully we design these assumptions, we must acknowledge that infill asymptotics inherently lead to coarser concentration bounds, with the “worst-case” infill asymptotics determining the precision of the Bernstein-type inequalities and the convergence rates of important nonparametric estimators. Unfortunately, this trade-off is an unavoidable consequence of the generality established by Assumptions \ref{Asp 1} and \ref{Asp 2}. Ultimately, the constraints imposed on sampled locations do not place any specific restrictions on their exact distribution or arrangement. Indeed, to the best of our knowledge, the assumptions outlined in Section 2.1 represent some of the mildest restrictions on DEI asymptotics found in the literature. At the same time, Proposition \ref{prop *} plays an important role in validating Theorem \ref{Theorem 1} because it provides  the required upper bounds for the variance of the “grouped variables” described in Section 2.1. These upper bounds serve as the foundation for applying Bernstein's blocking technique to irregularly spaced locations.}
\section{Bernstein-type Inequalities for NED}
{\par  This section introduces Bernstein-type inequalities for  irregularly spaced NED random fields. Compared to Hoeffding-type inequalities, a key advantage of Bernstein-type inequalities is their localization property, which stems from their reliance on the variance term of random variables. As a result, in asymptotic or non-asymptotic analysis of various estimators, Bernstein-type inequalities often facilitate the attainment of optimal convergence rates. Since the seminal works of \citeA{carbon1983inegalite} and \citeA{collomb1984proprietes} in the mid-1980s, there have been numerous extensions of Bernstein-type inequalities to various forms of dependent stochastic processes. Our paper specifically addresses NED random fields based on $\alpha$-mixing random fields. Accordingly, we focus here on reviewing the literature on Bernstein-type inequalities for $\alpha$-mixing processes. 
\par Building upon coupling and blocking techniques, \citeA{rio1995functional} and \citeA{liebscher1996strong} established Bernstein-type inequalities for stationary real-valued $\alpha$-mixing processes, assuming a broad range of dependence decay. Additionally, \citeA{bryc1996large} investigated the large deviation phenomenon for $\alpha$-mixing processes under exponential decay. Using a combination of Cantor set methodology and Bernstein's blocking technique, \citeA{merlevede2009bernstein} explored Bernstein-type inequalities under geometric $\alpha$-mixing conditions. More recently, efforts by \citeA{valenzuela2017bernstein} and \citeA{krebs2018bernstein} extended these findings to encompass geometric strong mixing spatial lattice processes and exponential graphs. However, all these studies have been limited to regularly spaced data, neglecting the exploration of data that is irregularly spaced. Surprisingly, there is a dearth of research into Bernstein-type inequalities or more comprehensive concentration inequalities for irregularly spaced data. \citeA{delyon2009exponential} proposed an exponential inequality for irregularly spaced mixing random fields using a martingale approach. However, the prerequisites for this result are challenging to specify and do not encompass the NED condition. As mentioned in Section 1, \citeA{xu2018sieve} presented the first direct result concerning exponential inequality for geometric irregularly spaced NED random fields, but their work focuses solely on pure domain-expand (DE) asymptotics.  
\subsection{Our Approach}
{\par Unlike \citeA{xu2018sieve}, who derived the exponential inequality by applying a covariance inequality along with various combinatorial techniques, our approach fully exploits the intrinsic properties of NED conditions, which naturally benefit NED processes. To elaborate, consider  Definition \ref{Def 1} and let $Z$ represent an irregularly spaced $L^p$ NED random field over an $\alpha$-mixing random field $\epsilon$. Our findings suggest that the probabilistic upper bound of the averaged partial sum process is governed by the concentration behavior exhibited by the random field $\epsilon$. This assertion aligns intuitively with Definition \ref{Def 1}. 
Importantly, for each $Z_{s_{iN}}\in Z$, the projection error -- i.e. the projection of $Z_{s_{iN}}$ onto the orthogonal complete space of $\mathcal{M}(\mathcal{F}_{iN}(r))$ -- converges to 0 in the $L^p$ norm as $r$ approaches $+\infty$. Here, $\mathcal{M}(\mathcal{F}_{iN}(r))$ denotes the collection of all real-valued $\mathcal{F}_{iN}(r)$-measurable functions. By allowing $r$ to increase monotonically with the sample size $N$, the random field (array) $\epsilon$ serves as an approximation of the field (array) $Z$ in the $L^{p}$ norm.}
\par Therefore, the assumption about the value of $p$ becomes crucial. Previous literature assumes $p\geq 2$, e.g., \citeA{lu2007local}, \citeA{li2012local} and \citeA{xu2018sieve}. More recently, \citeA{ren2020local} investigated the asymptotic properties of local linear quantile regression, albeit only under the condition $p=1$. In what follows, we now present our Bernstein-type inequality and its corollaries under $L^{p}$-NED conditions for any $p\geq 1$. 
\begin{theorem}
    \label{Theorem *}
    According to Definition \ref{Def 1}, suppose that the real-valued array $Z:=\{Z_{s_{iN}}:s_{iN}\in\Gamma_N\}$ is  $L^{p}$-NED($p\geq 1$) on array $\epsilon:=\{\epsilon_{s_{iN}}:s_{iN}\in\Gamma_N\}$, where $p\geq 1$. For any $i$ and $N$, assume $||Z_{s_{iN}}||_{\infty}\leq A$, $E|Z_{s_{iN}}|\leq \sigma_{1}$, $E[Z_{s_{iN}}^{2}]\leq \sigma^{2}$ and $E[Z_{s_{iN}}]=0$, for some $A,\sigma_{1},\sigma>0$. $\epsilon$ obeys the $\alpha$-mixing conditions introduced in Definition \ref{Def 2}. Let $\phi(x,y)=(x+y)^{\tau},\ \tau>0$ (see Definition \ref{Def 2}). Then, based on Assumptions \ref{Asp 1} to \ref{Asp 4+}, for any $t>0$ and $0<p_{N}<\beta(N)/2$, we have,
    \begin{align}
    \label{eq *1}
            \mathbb{P}(|N^{-1}S_{N}|>t)\leq& 2^{d_2+1}\exp\left(-\frac{Nd_{0}^{d}t^{2}}{C_{3}(\mathbf{V}d_{0}^{d}+C_{2}(p_N)^{d_2}At)}\right)\notag \\
    &+2^{d_2+\tau+\frac{1}{2}}11\left(\frac{2^{d_2+1}A}{t}+1\right)^{1/2}M(Q,N,\tau)\psi_{\epsilon}(\frac{p_{N}}{3})\notag \\
    &+2^{d_2}\left(\frac{2^{d_2+2}A(N)\psi_{Z}(\frac{p_{N}}{3})}{t}\right)^{p},
    \end{align}
    where $\mathbf{V}$ and $A(N)$ are introduced in Proposition \ref{prop *} and
    \begin{align*}
        &M(Q,N,\tau)=Q1[\tau=0]+(Q-1)^{1-\tau}N^{\tau}1[0<\tau<1]+N^{\tau}1[\tau\geq 1],\\
    &\frac{Nd_{0}^{d}}{2^{1.5d+d_{2}}H_{0}^{d_1}p_{N}^{d_2}}\leq Q=\frac{Vol(R_{N})}{H_{0}^{d_1}(2p_{N})^{d_2}}\leq \frac{KB_{N}}{2^{d_2}H_{0}^{d_1}p_{N}^{d_2}}.
    \end{align*}
   Here, $S_{N}=\sum_{i=1}^{N}Z_{s_{iN}}$ and $C_{2}$, $C_{3}$ are some universal constants independent of sample size $N$ and will be specified in the proof. Furthermore, if there exists some universal constant $C>0$ such that $\sigma_{1}\leq \frac{Ct}{2^{d_2+1}}$, we have
  \begin{align}
  \label{eq *2}
       \mathbb{P}(|N^{-1}S_{N}|>t)\leq& 2^{d_2+1}\exp\left(-\frac{Nd_{0}^{d}t^{2}}{C_{3}(\mathbf{V}d_{0}^{d}+C_{2}(p_N)^{d_2}At)}\right)\notag \\
    &+2^{d_2+1}11C^{1/3}M'(Q,N,\tau)\psi_{\epsilon}^{\frac{2}{3}}(\frac{p_{N}}{3})\notag\\
    &+2^{d_2}\left(\frac{2^{d_2+2}A(N)\psi_{Z}(\frac{p_{N}}{3})}{t}\right)^{p},
  \end{align}
  where $M'(Q,N,\tau)=Q1[\tau=0]+(Q-1)^{1-\frac{2\tau}{3}}N^{\frac{2\tau}{3}}1[0<\tau<\frac{3}{2}]+N^{\tau}1[\tau\geq \frac{3}{2}]$.
\end{theorem}
\par Unlike the commonly defined uniform NED condition (see Definition 1, \citeA{jenish2012spatial}), the assumption $\max_{i,N}\alpha_{s_{iN}}\leq A(N)$ is more general and may be more useful in certain applications. For example, suppose $Z=\{Z_{s_{iN}}\}_{s_{iN}\in\Gamma_N}$ is a $\mathcal{Z}$-valued $L^p$ irregularly spaced NED random field on $\epsilon=\{\epsilon_{s_{iN}}\}_{s_{iN}\in\Gamma_N}$, $\mathcal{Z}\subset\mathbf{R}^{D}$. Let $f:\mathcal{Z}\rightarrow \mathbf{R}$ be a Lipschitz continuous function with module $Lip(f)$. According to Proposition \ref{prop 4}, we know that $\{f(Z_{s_{iN}})\}_{s_{iN}\in\Gamma_N}$ can be regarded as a real-valued NED process on $\epsilon$ and $||f(Z_{s_{iN}})-E[f(Z_{s_{iN}})|\mathcal{F}_{iN}(r)]||_{p}\leq 2Lip(f)\psi_{Z}(r)$. However, $Lip(f)$ may sometimes diverge to infinity as $N\rightarrow +\infty$. A typical example is a Lipschitz continuous kernel function, $K(\frac{Z_{s_{in}}-z}{h})$, whose Lipschitz module is $O(1/h)$. 
\par Note that our Theorem \ref{Theorem *} applies to any given $d_{0} > 0$, as introduced in Assumption \ref{Asp 2}, meaning that it holds true for both DE and DEI asymptotics. However, a distinction arises when considering the impact of infill asymptotics and the effective dimension of $\Gamma_{N}$. In contrast to previous work on Bernstein-type inequalities for weakly dependent data (e.g., \citeA{rio1995functional}, \citeA{liebscher1996strong}, and \citeA{yao2003exponential}), where sharpness primarily depends on terms such as $t$, block size ($p_{N}$ in our case), variance term ($\mathbf{V}$ in our case) and the dimension of index set ($d$ in our paper), our inequality exhibits a key difference. The block size $p_{N}$ only affects the sharpness at the power of $d_{2}$ instead of $d$ (e.g., \citeA{yao2003exponential}). Another relevant point is that $p_N$ in Theorem \ref{Theorem *} serves as a tuning parameter within the range $(0, \beta(N))$, where $\beta(N)$ serves as the lower bound of the growth condition of $\frac{1}{2}\min_{1\leq k\leq d_{2}}H_{k}^{N}$. However, as noted earlier, the restriction on $\beta(N)$ is not fundamentally necessary because it is introduced to simplify the application of the blocking techniques used in our proof, particularly the discussion of how many irregularly spaced locations can be contained within a single block.
\par For the geometrical NED condition, we can easily obtain the following corollary which is a direct consequence of Theorem \ref{Theorem *}.
\begin{corollary}
    \label{Theorem 1}
    Based on the conditions and notation introduced in Theorem \ref{Theorem *} and $\mathbf{V1}$ defined in Proposition \ref{prop *}, let $Z$ be geometric $L^{p}$-NED on $\epsilon$ of order $(\nu, b,\gamma), \nu>1,$ and $\ b,\gamma>0$. Then, for any $t>0$ and some user defined parameter $\Theta>0$, if $\beta(N)$ introduced in Assumption \ref{Asp 4+} satisfies $\beta(N)\geq (\frac{6\Theta\log_{\nu}N}{b})^{\frac{1}{\gamma}}$, we have
\begin{align}
\label{eq:4}
     \mathbb{P}(|N^{-1}S_{N}|>t)\leq& 2^{d_2+1}\exp\left(-\frac{Nd_{0}^{d}t^{2}}{C_{3}(\mathbf{V1}d_{0}^{d}+C_{2}(\frac{3\Theta}{b}\log_{\nu}N)^{\frac{d_2}{\gamma}}At)}\right)\notag\\
    &+2^{d_2+\tau+\frac{1}{2}}11\left(\frac{2^{d_2+1}A}{t}+1\right)^{1/2}M(Q,N,\tau)N^{-\Theta}\notag \\
    &+2^{d_2}\left(\frac{2^{d_2+2}A(N)}{N^{\Theta}t}\right)^{p}.
\end{align}     
Similarly, when $\sigma_{1}\leq \frac{Ct}{2^{d_2+1}}$, by setting $\beta(N)\geq (\frac{9\Theta\log_{\nu}N}{b})^{\frac{1}{\gamma}}$, we have
\begin{align}
\label{eq:4'}
\mathbb{P}(|N^{-1}S_{N}|>t)\leq& 2^{d_2+1}\exp\left(-\frac{Nd_{0}^{d}t^{2}}{C_{3}(\mathbf{V1}d_{0}^{d}+C_{2}(\frac{9\Theta}{b}\log_{\nu}N)^{\frac{d_2}{\gamma}}At)}\right)\notag\\
    &+2^{d_2+1}11C^{1/3}M'(Q,N,\tau)N^{-\Theta}\notag \\
    &+2^{d_2}\left(\frac{2^{d_2+2}A(N)}{N^{\Theta}t}\right)^{p}.
\end{align}
\end{corollary}
The parameter $\Theta$ introduced here is user-defined. For instance, in the application of Corollary \ref{Theorem 1} to examine the uniform convergence rates of kernel-based estimators, the specific value of $\Theta$ is often linked to the dependence conditions and the complexity inherent in the underlying kernel classes. Consequently, assigning a higher value to $\Theta$ imposes more stringent requirements on the growth conditions of $\min_{1\leq k\leq d_2}H_{k}^{N}$.
Naturally, assuming that $\min_{1\leq k\leq d_{2}}H_{k}^{N}=O(N^{\theta})$, where $\theta\in (0,\frac{1}{d_2})$, allows Corollary \ref{Theorem 1} to hold for any $\Theta>0$.
\par Recall that Theorem A.1 of \citeA{xu2018sieve} is  one of the major motivations for this paper. To facilitate comparison with their work, we first summarize their result using our notation. Under Assumptions \ref{Asp 1} and \ref{Asp 2}, suppose the array $\{Z_{s_{iN}}:s_{iN}\in\Gamma_{N}\}$ has zero mean and uniform infinite norm $A>0$. $U$ and $V$ are two disjoint finite subsets of $\Gamma_{N}$. Then, based on some upper bound of $Cov(\prod_{s_{iN}\in U} Z_{s_{iN}},\prod_{s_{jN}\in V}Z_{s_{jN}})$, which is of order $e^{|U|+|V|}\exp(-\rho(U,V))A^{|U|+|V|}$, they obtain 
\begin{align}
\label{eq:7}\mathrm{P}\left(\left|N^{-1}\sum_{i=1}^{N}Z_{s_{iN}}N\right|\geq t\right)\lesssim C_{1}^{*}\exp\left(-C_{2}^{*}N^{\frac{1}{2d+2}}t^2\right),
\end{align}
where both $C_{1}^{*}$ and $C_{2}^{*}$ are positive constants but $C_{2}^{*}$ are not related to variance information.
\par The advantages of (\ref{eq:4}) and (\ref{eq:4'}) are evident. First, they can be adapted to the counterexamples mentioned in the third limitation of (\ref{eq:7}) in Section 1 because only the effective dimension affects the sharpness of inequality up to a $\log_{\nu}N$ level, which is clearly sharper than (\ref{eq:7}). Second, Corollary \ref{Theorem 1} holds for the NED condition under the $L^{p}$-norm for any $p\geq 1$, whereas \eqref{eq:7} is only valid for $p\geq 2$.  Moreover, as shown in the applications,  (\ref{eq:4}) and (\ref{eq:4'}) can be readily applied to investigate the optimality of many widely used estimators. Furthermore, compared with (\ref{eq:7}), our results yield sufficiently sharp convergence rates for estimators whose variances become arbitrarily small asymptotically, such as the kernel density estimator. Lastly, (\ref{eq:4}) and (\ref{eq:4'}) are valid for any $\nu>0$, whereas according to the proof of (\ref{eq:7}), $\nu =e$ is a necessary condition. Importantly, Corollary \ref{Theorem 1} accommodates both DE and DEI asymptotics, improving the real world applicability and flexibility of our results. Overall, Corollary \ref{Theorem 1} of Theorem \ref{Theorem *} addresses several of the shortcomings of the result presented by \citeA{xu2018sieve}, as discussed in Section 1. Unfortunately there is no “free lunch”. \citeA{xu2018sieve} derived an exponential inequality on an irregularly spaced and weakly dependent random field with zero mean and bounded support, containing both NED and $\alpha$-mixing random fields as special cases. 
\par However, an important criterion for assessing any trade-off is determining whether the additional restriction is sufficiently mild while providing substantial benefits. Here, we should point out that Assumption \ref{Asp 3} is essential for introducing the effective dimension of $\Gamma_{N}$ and thus cannot be omitted. As discussed in Section 2, Assumption \ref{Asp 4}  provides more flexibility in the selection of locations. By setting $B_N$ in Assumption \ref{Asp 4} equal to $N$, the asymptotics of $\Gamma_{N}$ can accommodate a mixture of DE and DEI asymptotics. Thus, under certain circumstances, our Assumption \ref{Asp 4} is more general than the prerequisites of \eqref{eq:7}. Regarding Assumption \ref{Asp 4+}, as previously mentioned, it is introduced only to simplify our proof. 
\par Beyond geometric NED conditions, Proposition 1 from \citeA{xu2015maximum}  indicates that the spatial autoregressive Tobit model (SART) may only satisfy the algebraic NED condition. Thus, it is important to establish a Bernstein-type inequality tailored to algebraically irregularly spaced NED random fields. Motivated by this, we introduce Corollary \ref{Theorem 2}, as follows:
\begin{corollary}
    \label{Theorem 2}
    Based on the conditions and notation introduced in Theorem \ref{Theorem *} and $\mathbf{V2}$ introduced in Proposition \ref{prop *}, when $Z$ is algebraic $L^{p}$-NED on $\epsilon$ of order $\gamma$, for some $\gamma>0$. Suppose there exists some $\theta\in (0,1)$ such that $\beta(N)$ introduced in Assumption \ref{Asp 4+} satisfies $\beta(N)\geq N^{\frac{\theta}{d_{2}}}$, we have
 \begin{align}
    \label{eq: 8}
            \mathbb{P}(|N^{-1}S_{N}|>t)\leq& 2^{d_2+1}\exp\left(-\frac{Nd_{0}^{d}t^{2}}{C_{3}(\mathbf{V2}d_{0}^{d}+C_{2}(p_N)^{d_2}At)}\right)\notag \\
    &+2^{d_2+\tau+\frac{1}{2}}11\left(\frac{2^{d_2+1}A}{t}+1\right)^{1/2}M(Q,N,\tau)\left(\frac{p_{N}}{3}\right)^{-\gamma}\notag \\
    &+2^{d_2}\left(\frac{2^{d_2+2}A(N)3^{\gamma}}{tp_{N}^{\gamma}}\right)^{p}.
    \end{align}
  When there exists some universal $C>0$ such that $\sigma_{1}\leq \frac{Ct}{2^{d_2+1}}$, under the same restriction of $\beta(N)$, we obtain
  \begin{align}
      \label{eq: 9}
       \mathbb{P}(|N^{-1}S_{N}|>t)\leq& 2^{d_2+1}\exp\left(-\frac{Nd_{0}^{d}t^{2}}{C_{3}(\mathbf{V2}d_{0}^{d}+C_{2}(p_N)^{d_2}At)}\right)\notag \\
    &+2^{d_2+1}11C^{1/3}M'(Q,N,\tau)\left(\frac{p_{N}}{3}\right)^{-\frac{2\gamma}{3}}\notag \\
    &+2^{d_2}\left(\frac{2^{d_2+2}A(N)3^{\gamma}}{tp_{N}^{\gamma}}\right)^{p}.
  \end{align}
\end{corollary}
\par In contrast to the geometric NED and mixing conditions, we impose a relatively lenient constraint on the growth condition of $\beta(N)$ in this context. This is a natural choice given that the decay of the dependence coefficients now follows a polynomial rate. We do not specify the specific value of $\theta$ here because it is closely tied to the design of $p_N$. As demonstrated later in Section 4, when applying \eqref{eq: 8} and \eqref{eq: 9} to investigate the uniform convergence rates of kernel regression, the design of $p_{N}$ is intimately linked with smoothness, bandwidth choice, moment conditions of the error terms, and dependence assumptions. However, as we have frequently emphasized, this is merely a technical assumption introduced to simplify our proof without any loss of generality.

\par  As a parallel result of our Theorem \ref{Theorem *}, the following Corollary can be regarded as a Bernstein-type inequality for an array of real-valued random variables satisfying $\alpha$-mixing conditions, provided that $\Gamma_{N}$ satisfies Assumptions \ref{Asp 1} to \ref{Asp 4+}. 
\begin{corollary} (\textbf{$\alpha$-mixing random field})
    \label{Theorem 3}
    Suppose $\{Z_{s_{iN}}:s_{iN}\in \Gamma_{N}\}$ is an array of real-valued random variables satisfying the $\alpha$-mixing conditions introduced in Definition \ref{Def 2}. Assume that $E[Z_{s_{iN}}]=0$, $\max_{i,N}||Z_{s_{iN}}||_{1}\leq \sigma_1$, $\max_{i,N}||Z_{s_{iN}}||_{2}\leq \sigma$ and let $\max_{i,N}||Z_{s_{iN}}||_{\infty}\leq A$, for some $\sigma_1, \sigma$ and $A>0$. Then, under Assumptions \ref{Asp 1} to \ref{Asp 4+}, for any given $t>0$, we obtain
    \begin{align}
        \label{eq: 10}
        \mathbb{P}(|\sum_{i=1}^{N}\epsilon_{s_{iN}}|>Nt) &\leq2^{d_2+1}\exp\left(\frac{Nd_{0}^{d}t^{2}}{2^{d_2+3}(\frac{Q}{N}d_{0}^{d}\Omega+\frac{1}{6}C_{2}p_{N}^{d_2}At)}\right)\notag\\
        &+2^{\tau+0.5+d_2}11(\frac{2^{d_2}A}{t}+1)^{0.5}M(Q,N,\tau)\psi_{\epsilon}(p_N).
    \end{align}
    Furthermore, if there exists a universal $C>0$ such that $t<\frac{Ct}{2^{d_2}}$, then we obtain
    \begin{align}
        \label{eq: 11}
        \mathbb{P}(|\sum_{i=1}^{N}\epsilon_{s_{iN}}|>Nt) &\leq2^{d_2+1}\exp\left(\frac{Nd_{0}^{d}t^{2}}{2^{d_2+3}(\frac{Q}{N}d_{0}^{d}\Omega+\frac{1}{6}C_{2}p_{N}^{d_2}At)}\right)\notag\\
        &+2^{d_2+1}11C^{1/3}M'(Q,N,\tau)\psi^{\frac{2}{3}}_{\epsilon}(p_N).
    \end{align}
Here, $C_{2}$, $M(Q,N,\tau)$, and $M'(Q,N,\tau)$ are all derived directly from Theorem \ref{Theorem *}. The term $\Omega$ is defined as $\sup_{{\Lambda(p_{N})}}Var(\sum_{s_{iN}\in \Lambda(p_N)}Z_{s_{iN}})$, where the supremum is taken over all $\Lambda(p_{N})\subset \mathbf{R}^{d}$. The notation $\Lambda(p_{N})$ denotes rectangles in which there are $d_2$ out of $d$ edges that have a length of $p_N$ and the remaining edges have a length of $H_0$. 
\end{corollary}
\par Compared to the work of \citeA{adamczak2008tail}, \citeA{modha1996minimum}, and \citeA{merlevede2011bernstein} and \citeA{liebscher1996strong}, our Corollary \ref{Theorem 3} exhibits either improved sharpness or broader applicability across a wider range of settings. Although the result of \citeA{delyon2009exponential} could match the sharpness of the iid setting, Corollary \ref{Theorem 3} is unique in that it explicitly accounts for infill asymptotics. Additionally, our approach does not require certain preconditions related to projections (i.e., conditional expectations) that are necessary for \citeA{delyon2009exponential}.
\par  We now briefly discuss the possibility of relaxing the moment conditions used in Theorem \ref{Theorem *}. Similar to Theorem A.2 in \citeA{xu2018sieve}, when the random variables $\{Z_{s_{iN}}:s_{iN}\in\Gamma_N\}$ exhibit sub-exponential tail behavior, we can obtain an inequality that remains as sharp as Theorem \ref{Theorem *} up to a logarithmic factor. To explore this, we introduce the following assumption on sub-exponential tails.
\begin{assumption} (\textbf{Sub-exponential tails})
    \label{assumption exponential}
    We assume that  there exist some universal positive constants $\kappa$, $\eta$ independent of $N$ such that, for every $i$ and $N$, we have
    \begin{align*}
        \mathbb{P}(|Z_{s_{iN}}|>t)\leq 2\exp\Big(-\frac{t}{2\eta}\Big),\ t>\frac{\kappa^2}{\eta}.
    \end{align*}
\end{assumption}
Using this assumption and applying the truncation method, we obtain the following inequality:
\begin{corollary}
    \label{corollary 4} (\textbf{Sub-exponential tails})
    By replacing the uniform boundedness assumption in Theorem \ref{Theorem *} with Assumption \ref{assumption exponential} while keeping other conditions and notation unchanged, we obtain the following results. For any given $\Theta'>0$ and positive integer $q$, when $N\geq N_0:=\min\{N:\frac{\Theta'}{2q\eta}\log N\geq \frac{\kappa^2}{\eta}\}$, the following holds for any $t>0$,
    \begin{align*}
        \mathbb{P}(|N^{-1}S_{N}|>t)&\leq 2^{d_2+1}\exp\left(-\frac{Nd_{0}^{d}t^{2}}{C_{3}(\mathbf{V}d_{0}^{d}+C_{2}(p_N)^{d_2}\frac{\Theta'}{2q\eta}(\log N)t)}\right)\notag \\
    &+2^{d_2+\tau+\frac{1}{2}}11\left(\frac{2^{d_2}\Theta'\log N}{q\eta t}+1\right)^{1/2}M(Q,N,\tau)\psi_{\epsilon}(\frac{p_{N}}{3})\notag \\
    &+2^{d_2}\left(\frac{2^{d_2+2}A(N)\psi_{Z}(\frac{p_{N}}{3})}{t}\right)^{p}+\frac{2\max_{i,N}||Z_{s_{iN}}||_{p}}{N^{\Theta'}t},
    \end{align*}
    where $p$ is also a positive integer such that $\frac{1}{p}+\frac{1}{q}=1$.
\end{corollary}
Similar to the $\Theta$ in Corollary \ref{Theorem 1}, $\Theta'$ here is also a user-defined constant. As we pointed out earlier, even though the sub-exponential tail is indeed more general than uniform boundedness, it slightly reduces the sharpness of our inequality.
\section{Applications of Our Theoretical Results to Kernel Conditional Mean Estimators}

\par In this section, we investigate the uniform convergence rate of the kernel estimator for the conditional mean function under algebraic NED conditions of order $\gamma$ and DE asymptotics. In particular, assuming the regression function belongs to the Hölder class, we examine the attainability of optimal uniform convergence rates for the local linear estimator when the input vector has expanding support. Our results build upon and extend previous contributions, including those of \citeA{hansen2008uniform}, \citeA{jenish2012nonparametric} and \citeA{chen2015optimal}.

\par For regularly spaced time series data, \citeA{hansen2008uniform} provided a thorough discussion of the uniform convergence rates of various kernel-based estimators under $\alpha$-mixing conditions and strict stationarity. \citeA{vogt2012nonparametric} later identified an ambiguous aspect of Hansen’s proof concerning the uniform convergence rate of a generic kernel estimator and delivered a revised proof for the uniform convergence rate of the Nadayara-Watson estimator using locally stationary and $\alpha$-mixing processes. \citeA{li2012local}  revisited this topic by extending mixing conditions to NED conditions. However, none of these studies explicitly address the sufficient conditions required to achieve optimal uniform convergence rates for kernel-based estimators.

\par For regularly spaced spatial data under $\alpha$-mixing conditions,  \citeA{hallin2004local} and \citeA{hallin2009local} showed convergence rates and asymptotic normality for local linear estimators of conditional mean and quantile functions. The problem becomes considerably more complex when dealing with irregularly spaced spatial data. To our knowledge, \citeA{jenish2012nonparametric} is the only study that discusses the uniform convergence rates of kernel-based estimators for the conditional mean function in the context of irregularly spaced spatial data with NED conditions. However, the uniform convergence rate they derive is non-standard and fails to achieve the optimal rate even with an appropriately chosen bandwidth. This gap in the literature motivates us to explore whether the optimal uniform convergence rate can be obtained for kernel smoothers estimating the conditional mean function under algebraic NED conditions and irregularly spaced locations. More specifically, we consider the following regression model:

\begin{align*}
    &Y_{s_{iN}}=m(X_{s_{iN}})+e_{s_{iN}},\\
    &E[e_{s_{iN}}|X_{s_{iN}}]=0,
\end{align*}
where $(X_{s_{iN}},Y_{s_{iN}})$ follow the same distribution as the random vector $(X,Y)\in\mathcal{X}\times\mathbf{R}$. 

\par As previously mentioned, we focus on the case in which $\mathcal{X}$ expands as sample size $N$ grows. Hence, we define $\mathcal{X}=\mathcal{B}_{D}(\rho_N)$, where $\mathcal{B}_{D}(\rho_N)$ denotes a $D$-dimensional Euclidean ball centered at the  origin with  radius $\rho_N$. We assume that $\rho_N=n^{\theta},$ for some $\theta>0$.
\par Given that $X_{s_{iN}}=x$, the local linear estimator of $m(x)$ is defined as:
\begin{align*}
    &\hat{m}(x)=(1,\textbf{0})\hat{\theta}(x),\ \hat{\theta}(x)=S_{N}(x)^{-1}T_{N}(x),\\
    &S_{N}(x)=\frac{1}{N}\sum_{i=1}^{N}K_{ih}(x)U_{i}(x)U_{i}(x)^{T},\ T_{N}(x)=\frac{1}{N}\sum_{i=1}^{N}K_{ih}(x)U_{i}(x)Y_{s_{iN}},\\
    &K_{ih}(x)=\frac{1}{h^{D}}\prod_{k=1}^{D}K\left(\frac{X^{k}_{s_{iN}}-x^{k}}{h}\right)=:\frac{1}{h^{D}}K\left(\frac{X_{s_{iN}}-x}{h}\right),\\
    &U_{i}(x)=\left(1,\left(\frac{X_{s_{iN}}-x}{h}\right)^{T}\right)^{T},
\end{align*}
where $\textbf{0}$ is a D-dimensional zero vector. Here, we assume that the product-type kernel satisfies Assumption \ref{Asp 22}. 
\begin{assumption}
\label{Asp 22}
  Suppose $K:\mathbf{R}\rightarrow\mathbf{R}$ is a function with compact support $[-1,1]$. For some fixed $\beta>0$, it satisfies
    \begin{itemize}
     \item [\textbf{K1.}] $\int K=1$, $||K||_{p}<\infty$ for any $1\leq p\leq \infty$.
    \item [\textbf{K2.}] $\int |t|^{\alpha}|K(t)|dt<\infty$. In the case of $[\alpha]\geq 1$, $\int t^{s}K(t)dt=0$ holds for any $0<s<\beta$.
    \item [\textbf{K3.}] The first-order derivative of $K$ is uniformly bounded on $[-1,1]$.
    \end{itemize}
\end{assumption}
Assumption \ref{Asp 22} coincides with the definition of a “$\beta$-valid” kernel as introduced by \citeA{tsybakov2004optimal}. Section 1.22 in \citeA{tsybakov2004introduction} shows how to construct such a kernel by using the orthonormal properties of Legendre polynomials. In the following, we outline important assumptions for our upcoming theorem.
\begin{assumption}(\textbf{Dependence and sampling})\\
\label{ASP 7}
Let the array $Z:=\{Z_{s_{iN}}:s_{iN}\in\Gamma_N\}$, where $Z_{s_{iN}}=(X_{s_{iN}},Y_{s_{iN}})$, be an algebraic $L^{p}$-NED ($p\geq 1$) process on the array $\epsilon:=\{\epsilon_{s_{iN}}:s_{iN}\in\Gamma_N\}$ of order $\gamma$. Here $\epsilon$ is an $\alpha$-mixing random field as defined in Definition \ref{Def 2}.
Assume that $\max_{iN}\alpha_{s_{iN}}=1$ and that $\psi_{Z}(t)=t^{-\gamma}$ and $\psi_{\epsilon}(t)=t^{-\gamma}$, for some $\gamma>0$. Additionally, assume $\tau=0$. i.e., $\phi(|U|,|V|)=1$. The set $\Gamma_{N}$ satisfies Assumptions \ref{Asp 1}-\ref{Asp 4+} with $\beta(N)\geq 2(N^{\frac{1}{p_0}}(\log N/N)^{\frac{\alpha}{2\alpha+D}})^{-\frac{1}{d_2}}$ with $d_0=1$. Here $p_0$ and $\alpha$ are the conditional moment conditions and smoothness index, respectively as  introduced in Assumptions \ref{ASP 8} and \ref{ASP 9}.
\end{assumption} 
\begin{assumption}(\textbf{Conditions of moments and density})
\label{ASP 8}
 Assume that, for any given $x\in\mathbf{R}^{D}$, $E[Y^{2}|X=x]\leq \sigma_{Y}^{2}$, $E[Y^{p_0}|X=x]< +\infty$, for some $p_0 >2$ and $\forall\ x\in\mathbf{R}^{D}$. Additionally assume that, for any $i$ and $j$, $X_{s_{iN}}$ and $(X_{s_{iN}},X_{s_{jN}})$ have Lebesgue densities, denoted as $f_{X_i}(u)$ and $f_{(X_i,X_j)}(x,y)$, respectively. Furthermore, we assume that the following conditions hold: 
\begin{enumerate}
    \item [M1.] $\max_{i,N}\sup_{x}f_{X_i}(x)<+\infty$.
    \item [M2.] For some $p_0>2$, $\max_{i,N}\sup_{x}E[|Y_{i}|^{p_0}|X_i=x]f_{X_i}(x)<+\infty$.
    \item [M3.] $\max_{i\neq j,N}\sup_{x}E[|Y_{i}Y_{j}||X_i=x,X_j=y]f_{(X_i,X_j)}(x,y).<+\infty$.
\end{enumerate}
\end{assumption}
\begin{assumption}(\textbf{Smoothness})\\
\label{ASP 9}
We assume that the conditional mean function $m(x)$ belongs to the Hölder space $\Sigma(\alpha, 1)$, for some $\alpha>0$. Here the definition of Hölder space coincides with Definition 1.2 in \citeA{tsybakov2004introduction}.  
\end{assumption}
\begin{assumption}(\textbf{Restrictions on $\gamma$, $p_0$ and expanding support})\\
    \label{ASP 10}
    We assume the following inequality holds for $\gamma$, $p_0$ and expanding support:
    \begin{enumerate}
        \item [R1.] $p_0>2+\frac{D}{\alpha}$.
        \item [R2.] $\gamma>\left\{\frac{((2\alpha+3)D+\alpha)p_0+(2\alpha+D)(1+2\theta Dp_0)}{(\alpha p_0-(2\alpha+D))}\right\}d_2$.
        \item [R3.] $\gamma>\frac{[(2\alpha+D)(\theta Dp_0/p+1)+p_0((\alpha+1)D/p+\alpha+D)]}{(\alpha p_0-(2\alpha+D))}d_2$
        \item [R4.] $\gamma\geq \{(\frac{1}{D}+1.5)\lor\frac{(2p_0-2)}{p_0-2}\}d_2$. 
    \end{enumerate}   
    Here $\theta$ is derived from the growth condition of $\rho_N$ introduced in Section 4, while $\alpha$ denotes the smoothness index from Assumption \ref{ASP 9}. The parameters $p$ and $p_0$ are based on the $L^p$-NED condition and M2 of Assumption \ref{ASP 8}. The term $d_2$ denotes the effective dimension of index set $\Gamma_N$.
\end{assumption}
For Assumption \ref{ASP 7}, setting $d_0=1$ is a standard assumption that does not lead to any loss of generality for DE asymptotics. To simplify the restriction on $\gamma$, we set $\tau=0$, making Assumption \ref{ASP 7} more directly comparable to the previously mentioned literature on time series data. However, extending the analysis to the case where $\tau>0$ can be achieved by following our proof. As discussed earlier, in Section 3, the restriction on $\beta(N)$ provides strong practical advantages when applying Corollary \ref{Theorem 2}, which is based on Bernstein's blocking technique under algebraic NED conditions. Assumption \ref{ASP 8} is directly borrowed from \citeA{hansen2008uniform} and \citeA{vogt2012nonparametric}. Assumption \ref{ASP 9} is also standard in the research on nonparametric regression.
\par Our discussion focuses on Assumption \ref{ASP 10} because it provides a sufficient condition for achieving the optimal uniform convergence rate. The parameter $\theta$ in R2 and R3 is derived from the definition of $\mathcal{B}_{D}(\rho_N)$ introduced on page 24 and serves as a user-defined tuning parameter. Its value is determined by certain moment conditions of the input vector introduced in (21) of Theorem 4 in \citeA{hansen2008uniform}. The requirement on $p_0$ ensures that, conditional on $X$, the order of the finite moment of output $Y$ must be strictly greater than $D/\alpha + 2$. Although sufficiently strong smoothness could allow this lower bound to be very close to $2$, this requirement is slightly more restrictive than that of Theorem 2.1 in \citeA{chen2015optimal}, in which $p_0$ could be equal to $D/\alpha+2$. However, because our proof strategy is based on Theorem 4.1 of \citeA{vogt2012nonparametric}, it is important for both of our results that this inequality holds strictly. A necessary condition for proving  Theorem 4.1 in  \citeA{vogt2012nonparametric} is that the term $S_T$ there diverges to $\infty$. When the bandwidth $h$ is set as the standard optimal bandwidth for uniform consistency (i.e., based on our notation introduced in Section 4.2, $h=(\log N/N)^{\frac{1}{2\alpha+D}}$), Vogt's truncation construction and block size imply that if $p_0=2+\frac{D}{\alpha}$, then $S_T\searrow 0$, which would invalidate the proof. Because our approach follows \citeauthor{vogt2012nonparametric}'s in this regard, we encounter a similar obstacle. Regarding R4 in Assumption \ref{ASP 10}, when $d_2=1$, the term $\gamma\geq \frac{(2p_0-2)}{p_0-2}$ coincides with the requirement for $\gamma$ in Theorem 1 of \citeA{hansen2008uniform}. Furthermore, R3 shows that imposing a stronger $L^p$-NED condition facilitates the relaxation of dependence assumptions. More importantly, R2, R3 and R4 together underscore the importance of introducing an effective dimension, represented by  $d_2$. For example, if our spatial data is drawn from a three dimensional random field ($d=3$) but, due to certain constraints, the locations are restricted  to a finite blow-up of a curve, 
such as a  mountain path, then $d_2$ in Assumption \ref{ASP 10} should be equal to $1$ rather than $3$. Recognizing this distinction considerably improves the generality of Assumption \ref{ASP 10}. Now we state our theorem as follows:
\begin{theorem}
\label{Theorem 4.2}
Under Assumptions \ref{Asp 22} to \ref{ASP 10}, by letting $h=(\log N/N)^{\frac{1}{2\alpha+D}}$, we have
\begin{align*}
    \sup_{x\in\mathcal{B}_{D}(\rho_N)}|\hat{m}(x)-m(x)|=O_{\mathbb{P}}\left(\left(\frac{\log N}{N}\right)^{\frac{\alpha}{2\alpha+D}}\right).
\end{align*} 
\end{theorem}

This result provides a set of sufficient conditions for achieving the optimal uniform convergence rate of the local linear estimator under algebraic NED conditions and DE asymptotics with irregularly spaced locations. As discussed earlier, our moment conditions and restrictions on $\gamma$ are more stringent than those of Theorem 2.1 of \citeA{chen2015optimal}. However, their Assumptions 2(i)(ii)(iii) in Theorem 2.1 require that the error terms are stationary martingale differences, which differs from our setting. However, as pointed out by \citeauthor{chen2015optimal}, martingale difference assumptions nest i.i.d. sequences while allowing for optimal uniform rates in sieve estimators under weak dependence conditions on regressors (inputs). This raises the question of whether a similarly refined approach could be adapted for kernel smoothers, which we aim to explore in future research. Another important difference is that \citeA{chen2015optimal} assume that the regressors satisfy $\beta$-mixing conditions, which are somewhat stricter than the $\alpha$-mixing conditions used in our paper. However, as noted in Theorem 4.2 of \citeA{chen2015optimal}, a key advantage of $\beta$-mixing is its connection to Berbee’s Lemma \cite{berbee1987convergence}. Compared to Bradley’s Lemma (Lemma 1.2 of \citeA{bosq2012nonparametric}) used in our proof, Berbee’s Lemma delivers a stronger bound on the error introduced by independent coupling. This suggests that, for the local linear smoother discussed in our paper, we may be able to relax some of our assumptions while still achieving optimal uniform convergence rates if we replace $\alpha$-mixing conditions with $\beta$-mixing.

\section{Discussion and Conclusion}

\par This paper examined two important topics under near-epoch dependent (NED) conditions: pure domain-expanding (DE) asymptotics and domain-expanding infill (DEI) asymptotics. Given that our definition of DEI asymptotics naturally encompasses DE asymptotics as a special case, we primarily focused on deriving Bernstein-type inequalities for irregularly spaced data under NED conditions and DEI asymptotics. As a corollary, we also derived a Bernstein-type inequality for irregularly spaced data under $\alpha$-mixing conditions and DEI asymptotics. A key aspect of our results is their adaptability to the effective dimension (see Definition 1) of the index set, which allows us to obtain much sharper upper bounds when the effective dimension is significantly smaller than the dimension of the index set.
\par Additionally, our inequalities were derived under a highly general formulation of DEI asymptotics because they depend only on the parameter $d_0$, which represents the strongest degree of infill asymptotics. Thus, our inequalities remain valid for any sample array $\Gamma_N$ satisfying Assumptions \ref{Asp 1} to \ref{Asp 4+}. In other words, by defining

$$\mathcal{S}_{N}=\{\Gamma_{N}:\Gamma_{N}\ \text{satisfies Assumptions \ref{Asp 1} to \ref{Asp 4+}}\},$$ 

our inequalities hold uniformly over family $S_N$. As a result, they naturally extend to the DEI asymptotics  introduced by \citeA{lu2014nonparametric}. Of course, sharper inequalities may be attainable under more specific restrictions on $\Gamma_N$, which presents an interesting avenue for future research. In terms of applications, we revisited the kernel-based conditional mean estimator under algebraic NED conditions and domain-expanding asymptotics. In particular, we identified a set of sufficient conditions for achieving the optimal uniform convergence rate of local linear smoothers for conditional mean function. Among these conditions, we highlighted the role of the effective dimension in determining the dependence conditions  necessary for attaining optimal rates.

\bibliographystyle{apacite}
\bibliography{bibliography.bib}
\newpage
\appendix

\section{Proof Related to Section 2}
\textbf{Proof of Proposition \ref{prop 1}}
\par For (i), we prove it by contradiction. Set $B(s,d_0/2)$ as an open ball in $\mathbf{R}^d$ whose center is $s\in\Gamma$ and radius is $d_0/2$. Obviously, $I(s,\sqrt{2}d_0/2)\subset B(s,d_0/2)$. If $B(s,d_0/2)$ contains two elements in $\Gamma_N$, say $s_{iN}$ and $s_{jN}$, a direct consequence is $||s_{jN}-s||\lor||s_{iN}-s||\leq d_0/2$. Thus, due to triangle inequality, $||s_{iN}-s_{jN}||\leq||s_{jN}-s||+||s_{iN}-s||\leq d_0$, which violates Assumption \ref{Asp 2}. Since $I(s,\sqrt{2}d_0/2)\subset B(s,d_0/2)$, we finish the proof. \\
For (ii), note that for $\forall h>0$, $\sqrt{2}h/d_{0}\leq[\sqrt{2}h/d_{0}]+1$. Thus $I(s,h)\subset I(s,([\sqrt{2}h/d_{0}]+1)\sqrt{2}d_{0}/2)$. Clearly, $I(s,([\sqrt{2}h/d_{0}]+1)\sqrt{2}d_{0}/2)$ can be partitioned into $([\sqrt{2}h/d_{0}]+1)^d$ cubes whose length of each is $\sqrt{2}d_{0}/2$. Then in light of Proposition 1 (i), $|I(s,h)\cap \Gamma_N|\leq ([\sqrt{2}h/d_{0}]+1)^d\leq (2\sqrt{2}h/d_0)^{d}.$\\
\par \textbf{Proof of Proposition \ref{prop 4}}
\par Define $g^{r}=E[g(Z_{s_{1N}},..,Z_{s_{|U|N}})|\mathcal{F}_{U}(r)]$ and $h^{r}=E[h(Z_{t_{1N}},..,Z_{t_{|V|N}})|\mathcal{F}_{V}(r)]$, where
$\mathcal{F}_{U}(r)=\sigma(\bigcup_{i=1}^{|U|}\mathcal{F}_{iN}(r))$, $\mathcal{F}_{V}(r)=\sigma(\bigcup_{j=1}^{|V|}\mathcal{F}_{jN}(r))$. Then, by denoting $\triangle g=g-g^r$, $\triangle h=h-h^r$, we have
\begin{align*}
    &Cov(g(Z_{s_{1N}},...,Z_{s_{|U|N}}),h(Z_{t_{1N}},...,Z_{t_{|V|N}}))\\
    &:= Cov(f,g)\notag = Cov(\triangle g+g^r, \triangle h+h^r)\\
    &= Cov(\triangle g, \triangle h)+Cov(g^r, \triangle h)+Cov(\triangle g, h^r)+Cov(g^r, h^r).
\end{align*}
Note $Cov(\triangle g, \triangle h)=E(\triangle g\triangle h)$, since $E\triangle h=E\triangle g=0$. Meanwhile,
\begin{align*}
    Cov(\triangle g, h^r)&= E[(\triangle g-E\triangle g)(h^r-Eh^r)]=E[\triangle g(h^r-Eh^r)]\\
    &=E[\triangle g(h^r+\triangle h-\triangle h+E\triangle h-E\triangle h-Eh^r)]\\
    &=E[\triangle g(h-Eh)]-E[\triangle g\triangle h]
\end{align*}
Similarly, $Cov(g^r, \triangle h)=E[(g-Eg)\triangle h]-E[\triangle g\triangle h]$, \\
which implies $Cov(g,h)=E[\triangle h(g-Eg)]+E[\triangle g(h-Eh)]-E[\triangle g\triangle h]+Cov(g^r,h^r)$.\\
Therefore, for $\forall p\geq 2$ and $q=\frac{p}{p-1}$,
\begin{align*}
    |Cov(g,h)|&\leq E|\triangle g\triangle h|+E|\triangle h(g-Eg)|+E|\triangle g(h-Eh)|+|Cov(g^r,h^r)|\\
              &\leq ||\triangle g||_{q}||\triangle h||_{p}+||\triangle h||_q||g-Eg||_p+||\triangle g||_q||h-Eh||_p+|Cov(g^r,h^r)|\\
              &:= I_1+I_2+I_3+I_4
\end{align*}
For $I_4$, since $g^r$ and $h^r$ are $\mathcal{F}_{U}(r)$ and $\mathcal{F}_{V}(r)$ measurable, according to \citeA{davydov1968convergence}, for $\forall \delta>0$, we have
\begin{align*}
    I_4\leq \phi(|U|,|V|)\psi_{\epsilon}^{\frac{p+\delta-2}{p+\delta}}(\rho(U,V))||g^r||_{p+\delta}||h^r||_{p+\delta}.
\end{align*}
Together with
\begin{align*}
    ||g^r||_{p+\delta}&=(E|E[g|\mathcal{F}_{U}(r)]|^{p+\delta})^{\frac{1}{p+\delta}}\leq (E[E[|g|^{p+\delta}|\mathcal{F}_{U}(r)]])^{\frac{1}{p+\delta}}=||g||_{p+\delta},
\end{align*}
we obtain $I_4\leq \phi(|U|,|V|)\psi_{\epsilon}^{\frac{p+\delta-2}{p+\delta}}(\rho(U,V))||g||_{p+\delta}||h||_{p+\delta}.$
\par Please remark $||g-Eg||_{p}\leq||g||_p+(\int|Eg|^{p}dP)^{1/p}\leq 2||g||_{p}$, where the second inequality is due to Jensen inequality. Similarly, we have $||h-Eh||_p\leq 2||h||_{p}$. Since $p\geq 2$, $q=p/(p-1)\leq p$. Thus $||\triangle g||_q\leq ||\triangle g||_p$ and $||\triangle h||_{q}\leq ||\triangle h||_p$. Moreover, to bound $||\triangle g||_p$, by applying Minkovski inequality, we have 
\begin{align*}
    ||\triangle g||_{p}\leq ||g-\Bar{g}||_p+||\Bar{g}-g^r||_p:=G_1+G_2,
\end{align*}
where $ \Bar{g}=g(E[Z_{s_{iN}}|\mathcal{F}_{iN}(r)],i=1,..,|U|)$. For $G_1$, since g is coordinate wise Lipschitz continuous on $\mathcal{Z}^{|U|}$, we have 
\begin{align*}
    G_1\leq Lip(g)\begin{matrix}\sum_{i=1}^{|U|}\end{matrix}||Z_{s_{iN}}-E[Z_{s_{iN}}|\mathcal{F}_{iN}(r)]||_p\leq Lip(g)(\begin{matrix}\sum_{i=1}^{|U|}\end{matrix}\alpha_{s_{iN}})\psi_{Z}(r).
\end{align*}
As for $G_2$, similarly to bounding $||A_2||_p$ in Proposition 3, regarding that 
\begin{align*}
  |\Bar{g}-g^r|&=|E[\Bar{g}-g^r|\mathcal{F}_{U}(r)]|\leq E[|\Bar{g}-g^r||\mathcal{F}_{U}(r)]\\
    &\leq Lip(g)\begin{matrix}\sum_{i=1}^{|U|}\end{matrix}E[|Z_{s_{iN}}-E[Z_{s_{iN}}|\mathcal{F}_{iN}(r)]||\mathcal{F}_{U}(r)],
\end{align*}
we have
\begin{align*}
    G_2&\leq Lip(g)||\begin{matrix}\sum_{i=1}^{|U|}\end{matrix}E[|Z_{s_{iN}}-E[Z_{s_{iN}}|\mathcal{F}_{iN}(r)]||\mathcal{F}_{U}(r)]||_p\\
    &\leq Lip(g)\begin{matrix}\sum_{i=1}^{|U|}\end{matrix}\left(E|E[|Z_{s_{iN}}-E[Z_{s_{iN}}|\mathcal{F}_{iN}(r)]||\mathcal{F}_{U}(r)]|^{p}\right)^{\frac{1}{p}}\\
    &\leq Lip(g)\begin{matrix}\sum_{i=1}^{|U|}\end{matrix}\left(E(E[|Z_{s_{iN}}-E[Z_{s_{iN}}|\mathcal{F}_{iN}(r)]|^{p}|\mathcal{F}_{U}(r)])\right)^{\frac{1}{p}}\\
    &=Lip(g)\begin{matrix}\sum_{i=1}^{|U|}\end{matrix}||Z_{s_{iN}}-E[Z_{s_{iN}}|\mathcal{F}_{iN}(r)]||_p\\
    &\leq Lip(g)(\begin{matrix}\sum_{i=1}^{|U|}\end{matrix}\alpha_{s_{iN}})\psi_{Z}(r).
\end{align*}
Therefore, we prove that $||\triangle g||_p\leq 2Lip(g)(\begin{matrix}\sum_{j=1}^{|V|}\end{matrix}\alpha_{jN})\psi_{Z}(r)$. Similarly, we can obtain $||\triangle h||_p\leq 2Lip(h)(\begin{matrix}\sum_{j=1}^{|V|}\end{matrix}\alpha_{jN})\psi_{Z}(r)$ as well. Together with all the results above, we obtain
\begin{align*}
    & I_1\leq ||\triangle g||_p||\triangle h||_p\leq 4Lip(g)Lip(h)(\begin{matrix}\sum_{i=1}^{|U|}\end{matrix}\alpha_{s_{iN}})(\begin{matrix}\sum_{j=1}^{|V|}\end{matrix}\alpha_{jN})\psi_{Z}^{2}(r),\\
    & I_2\leq ||\triangle g||_p||h-Eh||_p\leq2||\triangle g||_p|| h||_p\leq 4Lip(g)(\begin{matrix}\sum_{i=1}^{|U|}\end{matrix})\psi_{Z}(r),\\
    & I_3\leq  4Lip(h)(\begin{matrix}\sum_{j=1}^{|V|}\end{matrix}\alpha_{jN})\psi_{Z}(r). \ \textbf{Q.E.D}
\end{align*}

\par \textbf{Proof of Proposition \ref{prop *}}
\par We note that Assumptions \ref{Asp 1} and \ref{Asp 2} are too general to offer us any specific information about the pattern of the spatial arrangement of locations. Hence, we need to follow the idea of transforming the “irregularly spaced” setting back to “regularly spaced” settings by constructing “grouped variables” and see whether we can control the moment conditions of these “grouped variables” according to Proposition \ref{prop 1}.
\par Based on the $d_0$ introduced in the Section 2 and Proposition \ref{prop 1}, we partition the set $\Lambda$ into many smaller half-open-half-closed cubes whose each edge has a length $d_0/\sqrt{2}$. Then, according to Proposition \ref{prop 1}, each small cube could only contain no more than $1$ location. We denote these cubes as $\{c_{l}:1\leq l\leq n^{*}\}$. Here $n^{*}=n_{0}^{d-d_2}\prod_{k=1}^{d_2}n_{k}$, where, for $k=0,1,2,\dots,d_2$, 
\begin{align*}
       n_{k}= \left\{\begin{array}{rcl}
        \frac{\sqrt{2}h_{k}}{d_0},  &   &  \frac{\sqrt{2}h_{k}}{d_0}\in \mathbf{N};\\
        \left[ \frac{\sqrt{2}h_{k}}{d_0}\right]+1,  &   & \text{ otherwise}. 
    \end{array}
    \right.    
\end{align*}
Obviously these $c_{l}$'s are regularly spaced. Thus, for each given $d_0$, we could use the multiplication between the quantity of these small cubes(always integer-valued) and $d_0$ to express the lower bound of the distance between any $s_{iN}\neq s_{jN}\in \Gamma_{N}$. 
\par According to some basic algebra of variance, we have
\begin{align*}
    Var(\sum_{s_{iN}\in\Lambda}Z_{s_{iN}})\leq N(\Lambda)\sigma^{2}+\sum_{s_{iN}\neq s_{jN}}|Cov(Z_{s_{iN}},Z_{s_{jN}})|.
\end{align*}
By denoting $Cov(Z_{s_{iN}},Z_{s_{jN}})$ as $C_{ij}$, together with Proposition \ref{prop 4} and the condition that $\max_{i,N}\alpha_{s_{iN}}\leq A(N)$, we have
\begin{align*}
   |C_{ij}|&\leq 4A(N)\psi_{Z}^{2}(\frac{||s_{jN}-s_{iN}||}{3})\\
   &+8A(N)\psi_{Z}(\frac{||s_{jN}-s_{iN}||}{3})\sigma+2^{\tau}\psi_{\epsilon}^{\frac{\delta}{2+\delta}}(\frac{||s_{jN}-s_{iN}||}{3})\sigma_{2+\delta}^{2}\\
    &:=C_{ij}^1+C_{ij}^{2}+C_{ij}^{3}.
\end{align*}
Define $\{\tau_{N}:N\in\mathbf{N}\}$ as an integer-valued number sequence such that $\lim_{N}\tau_{N}\nearrow +\infty$. Define 
\begin{align*}
    R_{iN}=\{s_{jN}=(s_{jN}^{1},\dots,s_{jN}^{d})\in\Lambda:|s_{jN}^{k}-s_{iN}^{k}|\leq \frac{\tau_{N}d_{0}}{\sqrt{2}},\ 1\leq k\leq d_2 \}
\end{align*}
and $R_{iN}^{C}=\Lambda-R_{iN}$. An obvious result is that, for each $i$ and $N$, $R_{iN}$ could contain at most $n_{0}^{d_1}(2\tau_{N})^{d_2}$ locations in $\Gamma_{N}$, where $d_1=d-d_2$.
Them we have
\begin{align*}
    \sum_{s_{jN}\neq s_{iN}}|C_{ij}|=\sum_{s_{iN}\in\Lambda}\left(\sum_{s_{jN}\in R_{iN}}|C_{ij}|+\sum_{s_{jN}\in R_{iN}^{C}}|C_{ij}|\right).
\end{align*}
Obviously, $\sum_{s_{jN}\in R_{iN}}|C_{ij}|\leq 2^{d_2}n_{0}^{d_1}\tau_{N}^{d_{2}}\Sigma$, where $\Sigma=\max_{i,j}(E|Z_{s_{iN}}Z_{s_{jN}}|+E|Z_{s_{iN}}|E|Z_{s_{jN}}|)$. Meanwhile, it's easy to obtain that
\begin{align*}
    \sum_{s_{jN}\in R_{iN}^{C}}|C_{ij}|\leq \sum_{l=1}^{3}(\sum_{s_{jN}\in R_{iN}^{C}}C_{ij}^{l}):=\sum_{l=1}^{3}C_{i}^{l}.
\end{align*}
To treat $C_{i}^{l}$'s in a unified way, we investigate 
\begin{align*}
    \sum_{s_{jN}\in R_{iN}^{C}}\psi^{s}(\frac{||s_{jN}-s_{iN}||}{3}),\ \text{for some $s>0$,} 
\end{align*}
where $\psi(r)$ is some positive-valued function satisfying $\lim_{r\rightarrow +\infty}\psi(r)\searrow 0$. Recall that we have already partitioned the whole $R_{iN}$ into $n^*$ smaller cubes and each cube contains at most one location. Obviously this arguments works the same to $R_{iN}^{C}$. This means that, for any $s_{jN}\in R_{iN}^{C}$, due to the fact that these smaller cubes are half-open-half-closed, there exists only one $c_{l}$ such that $s_{jN}\in c_{l}$ and the volume of this cell is $(d_{0}/\sqrt{2})^{d}$. Together with the definition of Rieman integration, we have the following results.
\begin{align*}
     \sum_{s_{jN}\in R_{iN}^{C}}\psi^{s}(\frac{||s_{jN}-s_{iN}||}{3})&= \sum_{s_{jN}\in R_{iN}^{C}}\psi^{s}(\frac{||s_{jN}-s_{iN}||}{3})\left(\frac{d_{0}/\sqrt{2}}{d_{0}/\sqrt{2}}\right)^{d}\\
     &=d_{0}^{-d}\sum_{s_{jN}\in R_{iN}^{C}}\psi^{s}(\frac{||s_{jN}-s_{iN}||}{3})d_{0}^{d}\\
     &\leq d_{0}^{-d}\int_{[0,h_{0}]^{d_1}}\int_{U_{iN}}\psi^{s}(\frac{||u-s_{iN}||}{3})du.
\end{align*}
where $U_{iN}=\prod_{k=1}^{d_2}(-\infty,s_{iN}^{k}-\frac{\tau_{N}d_{0}}{\sqrt{2}}]\cup[s_{iN}^{k}+\frac{\tau_{N}d_{0}}{\sqrt{2}}+\infty)\subset \mathbf{R}^{d_2}$. Please note $(\sum_{k=1}^{d_2}(s_{iN}^{k}-u^{k})^{2})^{1/2}\leq ||u-s_{iN}||$. According to the monotonicity of $\psi(\cdot)$, we have
\begin{align*}
    &\sum_{s_{jN}\in R_{iN}^{C}}\psi^{s}(\frac{||s_{jN}-s_{iN}||}{3})\leq d_{0}^{-d}\int_{[0,h_{0}]^{d_1}}\int_{U_{iN}}\psi^{s}(\frac{||u-s_{iN}||}{3})du\\
    &\leq d_{0}^{-d}h_{0}^{d_1}\int_{U_{iN}}\psi^{s}\left(\frac{\sqrt{\sum_{k=1}^{d_2}(u^{k}-s_{iN}^{k})^{2}}}{3}\right)du\leq d_{0}^{-d}h_{0}^{d_1}\int_{\{u\in\mathbf{R}^{d_2}:||u||\geq \frac{\tau_{N}d_{0}}{\sqrt{2}}\}}\psi^{s}(\frac{u}{3})du\\
    &\leq  d_{0}^{-d}h_{0}^{d_1}(2\pi)^{d_2}\int_{\frac{\tau_{N}d_{0}}{\sqrt{2}}}^{+\infty}\rho^{d_2-1}\psi^{s}(\frac{\rho}{3})d\rho,
\end{align*}
where the last inequality is due to polar coordinate transformation. Thus, we obtain
\begin{align*}
    &C_{i}^{1}\leq(2\pi)^{d_2}h_{0}^{d_1}d_{0}^{-d}4A(N)\int_{\frac{\tau_{N}d_{0}}{\sqrt{2}}}^{+\infty}\rho^{d_2-1}\psi_{Z}^{2}(\frac{\rho}{3})d\rho:=(2\pi)^{d_2}h_{0}^{d_1}d_{0}^{-d} \mathbf{T}_{1};\\
    &C_{i}^{2}\leq (2\pi)^{d_2}h_{0}^{d_1}d_{0}^{-d}8A(N)\sigma\int_{\frac{\tau_{N}d_{0}}{\sqrt{2}}}^{+\infty}\rho^{d_2-1}\psi_{Z}(\frac{\rho}{3})d\rho:=(2\pi)^{d_2}h_{0}^{d_1}d_{0}^{-d} \mathbf{T}_{2};\\
    &C_{i}^{3}\leq (2\pi)^{d_2}h_{0}^{d_1}d_{0}^{-d} 2^{\tau}\sigma_{2+\delta}^{2}\int_{\frac{\tau_{N}d_{0}}{\sqrt{2}}}^{+\infty}\rho^{d_2-1}\psi_{\epsilon}^{\frac{\delta}{2+\delta}}(\frac{\rho}{3})d\rho:=(2\pi)^{d_2}h_{0}^{d_1}d_{0}^{-d} \mathbf{T}_{3}.
\end{align*}
Above all, we have 
\begin{align*}
    \sum_{s_{iN}\neq s_{jN}}|C_{ij}|\leq 2^{d_2}N(\Lambda)(\pi^{d_2}h_{0}^{d_1}d_{0}^{-d} (\mathbf{T}_{1}+\mathbf{T}_{2}+\mathbf{T}_{3})+n_{0}^{d_1}\tau_{N}^{d_2}\Sigma).
\end{align*}
This yields
\begin{align*}
  Var(\sum_{s_{iN}\in\Lambda}Z_{s_{iN}})\leq N(\Lambda)(\sigma^{2}+2^{d_2}n_{0}^{d_1}\tau_{N}^{d_2}\Sigma+ (2\pi)^{d_2}h_{0}^{d_1}d_{0}^{-d} (\mathbf{T}_{1}+\mathbf{T}_{2}+\mathbf{T}_{3})).
\end{align*}
Now we try to bound terms $\mathbf{T}_{1}$ to $\mathbf{T}_{3}$ when $\psi_{Z}(r)=\psi_{\epsilon}(r)=\nu^{-br^{\gamma}}$ and $\psi_{Z}(r)=\psi_{\epsilon}(r)=r^{-\gamma}$ respectively, for some $\nu>1$, $b,\gamma>0$. 
\par When $\psi_{Z}(r)=\psi_{\epsilon}(r)=\nu^{-br^{\gamma}}$, we first investigating the following integral.
\begin{align*}
    I(s,a)=\int_{a}^{+\infty}\rho^{d_2-1}\nu^{-bs(\rho/3)^{\gamma}}d\rho,\ s>0, a\geq 1.
\end{align*}
Let $L^*=\min\{L\in\mathbf{N}:\frac{d_2}{\gamma}-L\leq 0\}$. Then, by using integration by parts, we can show that
\begin{align*}
    I(s,a)\leq \frac{1}{\gamma}&\left(\sum_{k=1}^{L^*-1}(\frac{3^{\gamma}}{bs})^{k}\left(\frac{d_2}{\gamma}(\frac{d_2}{\gamma}-1)\dots(\frac{d_2}{\gamma}-k+1)\right)a^{d_2-k\gamma}\nu^{-bs(a/3)^{\gamma}}\right.\\
    &\left.+ (\frac{3^{\gamma}}{bs})^{L^{*}}\left(\frac{d_2}{\gamma}(\frac{d_2}{\gamma}-1)\dots(\frac{d_2}{\gamma}-L^{*}+1)\right)\nu^{-bs(a/3)^{\gamma}}              
        \right).
\end{align*}
Thus, there exists some universal constant $C^{*}>0$ such that, for any $s>0$ and $a\geq 1$, we have
\begin{align*}
    I(s,a)\leq C^{*}a^{d_2-\gamma}\nu^{-bs(a/3)^{\gamma}}.
\end{align*}
Then, according to some simple algebra, we get upper bound $\mathbf{B1}$ defined in Proposition \ref{prop +}.
\par When $\psi_{Z}(r)=\psi_{\epsilon}(r)=r^{-\gamma}$, we investigate the following integration.
\begin{align*}
    I(s,a)&=\int_{a}^{+\infty}\rho^{d_2-1}\rho^{-s\gamma}d\rho,\ a\geq 1,\ s>0, \gamma>\frac{d_2}{s}.
\end{align*}
Based on some simple calculation, we get
\begin{align*}
    I(s,a)=(s\gamma-d_2)^{-1}a^{-(s\gamma-d_{2})}.
\end{align*}
Here $s=1,2,\frac{\delta}{2+\delta}$. Thus we require $\gamma>\frac{d_{2}(2+\delta)}{\delta}$. The result above directly leads us to upper bound $\mathbf{B2}$ described in Proposition \ref{prop *}.
\section{Proof Related to Section 3}
\begin{lemma}\rm{(\textbf{Lemma 1.2} of \citeA{bosq2012nonparametric})}
\label{lemma 3}
\textit{Let $(X,Y)$ be a $\mathbf{R}^{n}\times\mathbf{R}$-valued random vector such that $Y\in L^{p}(P)$ for some $p\in[1,+\infty]$. Let $c$ be a real number such that $||Y+c||_{p}>0$, and $\xi\in [0,||Y+c||_{p}]$. Then, there exists a random variable $Y^{*}$ such that\\
\begin{enumerate}
    \item $P_{Y^{*}}=P_{Y}$ and $Y^{*}$ is independent of $X$.
    \item 
        $P(|Y-Y^{*}|>\xi)\leq 11\left(\frac{||Y+c||_{p}}{\xi}\right)^{\frac{p}{2p+1}}[\alpha(\sigma(X),\sigma(Y))]^{\frac{2p}{2p+1}}$.
\end{enumerate}}
\end{lemma}
\par \textbf{Proof of Theorem \ref{Theorem *}}
\par \textbf{Step 1}(Construction of blocks)
\par Without loss of generality, we assume the $R_{N}$ introduced in Section 2 and Assumption \ref{Asp 4} is as follow,
\begin{align*}
    R_{N}=\prod_{k=1}^{d_2}(0,H_{k}^{N}] \times (0,H_{0}]^{d_1},\ d_{1}+d_{2}=d.
\end{align*}
We also assume $\min_{1\leq k\leq d_{2}}H_{k}^{N}\geq H_{0}$ for any $N$ and $\lim_{N}H_{k}^{N}\nearrow +\infty$ for any $1\leq k\leq d_{2}$. According to our two definitions of effective dimension, we know $R_{N}$ has effective dimension no larger than $d_{2}$.
\par Let $p_{N}$ be a positive number associated with $N$ such that $0<p_{N}\leq \frac{1}{2}\min_{1\leq k\leq d_{2}}H_{k}^{N}$. Let $q_{kN}=H_{k}^{N}/2p_{N}$. Here, similar to the proof of Theorem 1.3 in \citeA{bosq2012nonparametric}, we do not need to restrain the image region of $q_{kN}$ to be positive integer set. Thus, $R_{N}$ is partitioned into $2^{d_{2}}$ collections of rectangles. Within each collection, there are $\prod_{k=1}^{d_{2}}q_{kN}:=Q$ rectangles. The volume of each 
rectangle is $p_{N}^{d_2}H_{0}^{d_1}$. We name these rectangles as $A_{jm}$, where $1\leq j\leq Q$ and $1\leq m\leq 2^{d_2}$. This is to say, $A_{jm}$ is the $j$-th rectangle of the $m$-th collection. Furthermore, we denote the selected locations contained by $A_{jm}$ as $N_{jm}$. i.e. $|A_{jm}\cap \Gamma_{N}|=N_{jm}$. Before we move to the next step, in order to simplify our statement later, we here exhibit some important relationship among $N_{jm}$, $Q$ and $N$ based on Assumption \ref{Asp 4}.
\begin{itemize}
    \item [1] $Q=\frac{\prod_{k=1}^{d_{2}}H_{k}^{N}}{2^{d_{2}}p_{N}^{d_2}}=\frac{Vol(R_{N})}{2^{d_2}H_{0}^{d_1}p_{N}^{d_2}}\leq C_{K}\frac{B_{N}}{p_{N}}$, where $C_{K}=\frac{K}{2^{d_2}H_{0}^{d_1}}$ and $B_{N}$ is equal to $N$ or $Nd_{0}^{d}$.
    \item [2] For each $m$, $\sum_{j=1}^{Q}N_{jm}:=N_{m}\leq N$. We also have $\frac{\max_{m}N_{m}}{N}=O(1)$.
\end{itemize}
Proof of the two points above is quite trivial hence omitted here. 
\par \textbf{Step 2}(Construction of grouped variables)
\par Let 
$$ G_{jm}=\sum_{i=1}^{N}1[s_{iN}\in A_{jm}]Z_{s_{iN}}=\sum_{s_{iN}\in A_{jm}}Z_{s_{iN}}$$ 
and, for any $0<r<p_{N}/2$, $A_{jm}^{r+}$ be an open $r$-blow up of rectangle $A_{jm}$. i.e. $A_{jm}^{r+}=\{s\in\mathbf{R}^{d}:||s-A_{jm}||<r\}$, where $||s-A_{jm}||=\inf\{s'\in A_{jm}:||s-s'||\}$. Define $\sigma(A_{jm}^{r+})$ as the sigma algebra generated by $\epsilon_{s_{iN}}$'s whose $s_{iN}\in A_{jm}^{r+}$. Finally, we define $G_{jm}^{r}=E[G_{jm}|\sigma(A_{jm}^{r+})]$. Then, we can immediately obtain the following fundamental properties.
\begin{itemize}
    \item [(1)] $S_{N}=\sum_{i=1}^{N}Z_{s_{iN}}=\sum_{m=1}^{2^{d_2}}\sum_{j=1}^{Q}G_{jm}$;
    \item [(2)] For any $r>0$, $E[G_{jm}]=E[G_{jm}^{r}]=0$, $||G_{jm}||_{\infty}\lor||G_{jm}^{r}||_{\infty} \leq N_{jm}A$ and $\sqrt{Var(G_{jm}^{r})}=||G_{jm}^{r}||_{2}\leq ||G_{jm}||_{2}=\sqrt{Var(G_{jm})}\leq \sqrt{N_{jm}\mathbf{V}} $, where $\mathbf{V}$ is defined in Proposition \ref{prop *}.
\end{itemize}
The third inequality in (2) is due to the Jensen inequality for conditional expectation and the fourth inequality is a direct result of Proposition \ref{prop *}. 
\par Above all, by denoting the law as $\mathbb{P}$, we have, for any given $t>0$,
\begin{align*}
    &\mathbb{P}\left(|N^{-1}S_{N}|>t\right)= \mathbb{P}\left(\left|\sum_{m=1}^{2^{d_2}}\sum_{j=1}^{Q}G_{jm}\right|>Nt\right)\\ \leq&\mathbb{P}\left(\left|\sum_{m=1}^{2^{d_2}}\sum_{j=1}^{Q}(G_{jm}-G_{jm}^{r})\right|>\frac{Nt}{2}\right)+\mathbb{P}\left(\left|\sum_{m=1}^{2^{d_2}}\sum_{j=1}^{Q}G_{jm}^{r}\right|>\frac{Nt}{2}\right)\\
    \leq& \sum_{m=1}^{2^{d_2}}\mathbb{P}\left(\left|\sum_{j=1}^{Q}\sum_{s_{iN}\in A_{jm}}\left(Z_{s_{iN}}-E[Z_{s_{iN}}|\sigma(A_{jm}^{r+})]\right)\right|>N\eta\right)+\sum_{m=1}^{2^{d_2}}\mathbb{P}\left(\left|\sum_{j=1}^{Q}G_{jm}^{r}\right|>N\eta\right)\\
    :=& \sum_{m=1}^{2^{d_2}}(\mathbf{P}_{m1}+\mathbf{P}_{m2}),\ \text{where $\eta=t/2^{d_2+1}$.}
\end{align*}
\par \textbf{Step 3}($\mathbf{P}_1$)
\par According to Markov inequality, by denoting $|\Gamma_{N}\bigcap (\bigcup_{j=1}^{Q}A_{jm})|$ as $N_{m}$, we immediately have, for any $p\geq 1$,
\begin{align*}
    \mathbf{P}_{m1}&\leq \mathbb{P}\left(\sum_{j=1}^{Q}\sum_{s_{iN}\in A_{jm}}|Z_{s_{iN}}-E[Z_{s_{iN}}|\sigma(A_{jm}^{r+})]|>N\eta\right)\\
    &\leq \frac{||\sum_{j=1}^{Q}\sum_{s_{iN}\in A_{jm}}(Z_{s_{iN}}-E[Z_{s_{iN}}|\sigma(A_{jm}^{r+})])||_{p}^{p}}{(N\eta)^{p}}\\
    &\leq \left(\frac{\sum_{j=1}^{Q}\sum_{s_{iN}\in A_{jm}}||Z_{s_{iN}}-E[Z_{s_{iN}}|\sigma(A_{jm}^{r+})]||_{p}}{N\eta}\right)^{p}\\
    &\leq \left(\frac{N_m}{N}\right)^{p}\left(\frac{\max_{1\leq j\leq Q}\max_{s_{iN}\in A_{jm}}||Z_{s_{iN}}-E[Z_{s_{iN}}|\sigma(A_{jm}^{r+})]||_{p}}{\eta}\right)^{p}.
\end{align*}
Now we focus on bounding $||Z_{s_{iN}}-E[Z_{s_{iN}}|\sigma(A_{jm}^{r+})]||_{p}$ based on the definition of NED condition. Recall the definition of $\mathcal{F}_{iN}(r)$. It is obvious that $\mathcal{F}_{iN}(r)\subset \sigma(A_{jm}^{r+})$ provided that $s_{iN}\in A_{jm}$. Hence, when $s_{iN}\in A_{jm}$, $E[Z_{s_{iN}}|\mathcal{F}_{iN}(r)]$ is measurable with respect to $\sigma(A_{jm}^{r+})$, which yields
\begin{align*}
    E[Z_{s_{iN}}|\mathcal{F}_{iN}(r)]=E[ E[Z_{s_{iN}}|\mathcal{F}_{iN}(r)]|\sigma(A_{jm}^{r+})].
\end{align*}
Then, for any given $A_{jm}$ and any $s_{iN}\in A_{jm}$,
\begin{align*}
    &||Z_{s_{iN}}-E[Z_{s_{iN}}|\sigma(A_{jm}^{r+})]||_{p}\\
    \leq &||Z_{s_{iN}}-E[Z_{s_{iN}}|\mathcal{F}_{iN}(r))]||_{p}+||E[Z_{s_{iN}}|\mathcal{F}_{iN}(r))]-E[Z_{s_{iN}}|\sigma(A_{jm}^{r+})]||_{p}\\
    \leq &A(N)\psi_{Z}(r)+||E[E[Z_{s_{iN}}|\mathcal{F}_{iN}(r))]|\sigma(A_{jm}^{r+})]-E[Z_{s_{iN}}|\sigma(A_{jm}^{r+})]||_{p}\\
    \leq &2A(N)\psi_{Z}(r).
\end{align*}
Based on this, we immediately obtain that
\begin{align*}
    \mathbf{P}_{m1}\leq \left(\frac{2A(N)\psi_{Z}(r)}{\eta}\right)^{p}=\left(\frac{2^{d_2+2}A(N)\psi_{Z}(r)}{t}\right)^{p}.
\end{align*}
\par \textbf{Step 3}($\mathbf{P}_{m2}$)
\par Note that, for each given $m$, $G_{jm}^{r}$ is a real-valued measurable function of $\epsilon_{s_{jN}}$'s located in $A_{jm}^{r+}$. Furthermore, it's easy to notice that, for any $1\leq j\neq j'\leq Q$, the Euclidean distance between $A_{jm}^{r+}$ and $A_{j'm}^{r+}$ is no less than $p_{N}-2r$. Then, given an $n+1$-tuple $(j_{1}m, j_{2}m,\dots j_{n}m,j_{n+1}m)$, by denoting $|A_{jm}^{r+}\bigcap \Gamma_{N}|$ as $N_{jm}^{r}$, we have
\begin{align*}
    &\alpha(\sigma(G_{j_1,m}^{r},\dots,G_{j_n,m}^{r}), \sigma(G_{j_{n}m}^{r}))\\
    \leq& \alpha (\sigma(\bigcup_{k=1}^{n}A_{j_{k}m}^{r+}),\sigma(A_{j_{n+1}m}^{r+}))\leq (\sum_{k=1}^{n+1}N_{j_{k}m}^{r})^{\tau}\alpha(p_{N}-2r).
\end{align*}
Then, according to Lemma\ref{lemma 3}, we can construct an array of mutually independent random variables, denoted as $\{W_{jm}:1\leq j\leq Q, 1\leq m\leq 2^{d_2}\}$, such that $W_{jm}$ is identical to $G_{jm}^{r}$. This indicates that, for any given $m$ and $\xi_{jm}>0$, we have
\begin{align*}
    \mathbf{P}_{m2}\leq \mathbb{P}\left(\left|\sum_{j=1}^{Q}W_{jm}\right|>N\eta-\sum_{j=1}^{Q}\xi_{jm}\right)+\sum_{j=2}^{Q}\mathbb{P}\left(\left|W_{jm}-G_{jm}^{r}\right|>\xi_{jm}\right).
\end{align*}
To bound $\mathbb{P}(|W_{jm}-G_{jm}^{r}|>\xi_{jm})$, we apply Lemma \ref{lemma 3}. Please note that $p$ here is allowed to be any real number in $[1,+\infty]$, since $G_{jm}^{r}$ has bounded support for every given $N$. Hence we can understand the moment “$p$” in Lemma \ref{lemma 3} as a parameter. Meanwhile, note the restrictions for $\xi_{jm}$ is that $0<\xi_{jm}<||G_{jm}^{r}+c||_{p}$, where $c$ is a constant such that $||G_{jm}^{r}+c||_{p}>0$. Since $\xi_{jm}>0$, we can assume $N_{jm}^{r}\geq 1$ without loss of generality, otherwise we directly obtain $\mathbb{P}(|W_{jm}-G_{jm}^{r}|>\xi_{jm})=0$. 
\par When $t$ is any positive number, we choose $p=\infty$. Then, considering that $||G_{jm}^{r}+c||_{\infty}=||G_{jm}^{r}||_{\infty}+c$,$||G_{jm}^{r}+c||_{\infty}>0$ holds for any $c>0$. Hence, we design $c=N_{jm}^{r}\eta$ and $\xi_{jm}=\frac{1}{2}N_{jm}^{r}\eta$, which meets the restrictions. Then, for any $j,m$, we have
\begin{align*}
    &\mathbb{P}(|W_{jm}-G_{jm}^{r}|>\xi_{jm})\\
    \leq &11\left(\frac{2N_{jm}^{r}(A+\eta)}{N_{jm}^{r}\eta}\right)^{1/2}(N_{jm}^{r}+N_{j-1,m}^{r})^{\tau}\psi_{\epsilon}(p_{N}-2r)\\
    = & 11\sqrt{2}\left(\frac{A}{\eta}+1\right)^{1/2}(N_{jm}^{r}+N_{j-1,m}^{r})^{\tau}\psi_{\epsilon}(p_{N}-2r).
\end{align*}
Then, we obtain that
\begin{align*}
    &\sum_{j=2}^{Q}\mathbb{P}(|W_{jm}-G_{jm}^{r}|>\xi_{jm})\\
    \leq &11\sqrt{2}\left(\frac{A}{\eta}+1\right)^{1/2}(\sum_{j=2}^{Q}(N_{jm}^{r}+N_{j-1,m}^{r})^{\tau})\psi_{\epsilon}(p_{N}-2r)\\
    \leq & 2^{\tau+\frac{1}{2}} 11\left(\frac{A}{\eta}+1\right)^{1/2}M(Q,N,\tau)\psi_{\epsilon}(p_{N}-2r),
\end{align*}
where $M(Q,N,\tau)=Q1[\tau=0]+(Q-1)^{1-\tau}N^{\tau}1[0<\tau<1]+N^{\tau}1[\tau\geq 1]$.
\par Furthermore, when there exists some universal $C>1$ such that $\sigma_1\leq C\eta$, we select $p=1$. Then, $||G_{jm}^{r}+c||_{1}\geq |c-||G_{jm}^{r}||_{1}|$. Note that $||G_{jm}^{r}||_{1}\leq N_{jm}^{r}\sigma_{1}$. By letting $c=2N_{jm}^{r}\sigma_{1}$, we obtain that $||G_{jm}^{r}+c||_{1}\geq N_{jm}^{r}\sigma_{1}>0$. Then, by setting $\xi_{jm}=\frac{1}{2C}N_{jm}^{r}\sigma_{1}$, we obtain
\begin{align*}
&\sum_{j=2}^{Q}\mathbb{P}(|W_{jm}-G_{jm}^{r}|>\xi_{jm})\\
\leq & 11\psi_{\epsilon}^{\frac{2}{3}}(p_{N}-2r)\sum_{j=2}^{Q}\left(\frac{3N_{jm}^{r}\sigma_{1}}{\frac{1}{2C}N_{jm}^{r}\sigma_{1}}\right)^{1/3}(N_{jm}^{r}+N_{j-1,m}^{r})^{\frac{2\tau}{3}}\\
\leq & 22C^{1/3}M'(Q,N,\tau)\psi_{\epsilon}^{\frac{2}{3}}(p_{N}-2r),
\end{align*}
where $M'(Q,N,\tau)=Q1[\tau=0]+(Q-1)^{1-\frac{2\tau}{3}}N^{\frac{2\tau}{3}}1[0<\tau<\frac{3}{2}]+N^{\tau}1[\tau\geq \frac{3}{2}]$.
\par For both cases mentioned above, we always have $N\eta-\sum_{j=1}^{Q}\xi_{jm}\geq \frac{1}{2}N\eta$. Therefore, according to the Bernstein-inequality for independent random variables, for each $m$, we have
\begin{align*}
    &\mathbb{P}\left(\left|\sum_{j=1}^{Q}W_{jm}\right|>N\eta-\sum_{j=1}^{Q}\xi_{jm}\right)\leq \mathbb{P}\left(\left|\sum_{j=1}^{Q}W_{jm}\right|>\frac{1}{2}N\eta\right)\\
    \leq & 2\exp\left(-\frac{(N\eta)^{2}}{8(\sum_{j=1}^{Q}Var(G_{jm}^{r})+\frac{1}{3}(\max_{j}N_{jm}^{r})AN\eta)}\right).
\end{align*}
Recall that $Var(G_{jm}^{r})\leq Var(G_{jm})\leq N_{jm}\mathbf{V}$, where $\mathbf{V}$ is defined in Proposition \ref{prop *}. Then, $\sum_{j=1}^{Q}Var(G_{jm}^{r})\leq N\mathbf{V}$. Meanwhile, according to Proposition \ref{prop 1}, we can easily show $\max_{j,m}N_{jm}^{r}\leq (2\sqrt{2})^{d}H_{0}^{d_1}(p_{N}+r)^{d_2}/d_{0}^{d}$. By setting $r=\frac{p_{N}}{3}$, together with the fact that $\eta=t/2^{d_2+1}$, we have
\begin{align*}
&\mathbb{P}\left(\left|\sum_{j=1}^{Q}W_{jm}\right|>\frac{1}{2}N\eta\right)\leq 2\exp\left(-\frac{Nt^{2}}{C_{3}(\mathbf{V}+C_{2}(p_N)^{d_2}d_{0}^{-d}A\eta)}\right),\\
&\text{where}\ C_{2}=(\frac{4}{3})^{d_2}2^{d_1-1+\frac{d}{2}}H_{0}^{d_1},\ C_{3}=2^{2d_2+5}.
\end{align*}
Above all, we manage to show, for any given $t>0$, 
\begin{align*}
    \mathbb{P}(|N^{-1}S_{N}|>t)\leq& 2^{d_2+1}\exp\left(-\frac{Nt^{2}}{C_{3}(\mathbf{V}+C_{2}(p_N)^{d_2}d_{0}^{-d}A\eta)}\right)\\
    &+2^{d_2+\tau+\frac{1}{2}}11\left(\frac{2^{d_2+1}A}{t}+1\right)^{1/2}M(Q,N,\tau)\psi_{\epsilon}(\frac{p_{N}}{3})\\
    &+2^{d_2}\left(\frac{2^{d_2+2}A(N)\psi_{Z}(\frac{p_{N}}{3})}{t}\right)^{p}.
\end{align*}
More specifically, if there exists a universal $C>0$ such that $\sigma_{1}\leq C\eta=Ct/2^{d_2+1}$, we have
\begin{align*}
     \mathbb{P}(|N^{-1}S_{N}|>t)\leq& 2^{d_2+1}\exp\left(-\frac{Nt^{2}}{C_{3}(\mathbf{V}+C_{2}(p_N)^{d_2}d_{0}^{-d}A\eta)}\right)\\
    &+2^{d_2+1}11C^{1/3}M(Q,N,\tau)\psi_{\epsilon}^{\frac{2}{3}}(\frac{p_{N}}{3})\\
    &+2^{d_2}\left(\frac{2^{d_2+2}A(N)\psi_{Z}(\frac{p_{N}}{3})}{t}\right)^{p}.
\end{align*}\\

\par \textbf{Proof of Corollary \ref{Theorem 3}}
\par Since the proof is basically the same as the proof of Theorem \ref{Theorem *}, we only sketch the main steps here and also inherit the notation. \textbf{Step 1} is exactly the same as that of the proof of Theorem \ref{Theorem *}. In \textbf{Step 2}, we only need to construct random variables $G_{jm}$. $G_{jm}^{r}$ is no longer needed since we no longer need to deal with the projection errors caused by the definition of NED conditions. Then, based some basic properties of measure, we can show
\begin{align*}
    &\mathbb{P}(|N^{-1}\sum_{i=1}^{N}\epsilon_{s_{iN}}|>t)=\mathbb{P}\left(\left|\sum_{m=1}^{2^{d_2}}\sum_{j=1}^{Q}G_{jm}\right|>Nt\right)\\
    \leq&  \sum_{m=1}^{2^{d_2}}\mathbb{P}\left(\left|\sum_{j=1}^{Q}G_{jm}\right|>N\zeta\right):=\sum_{m=1}^{2^{d_2}}P_{m},
\end{align*}
where $\zeta=t/2^{d_2}$. Similar to \textbf{Step 3} in the proof of Theorem \ref{Theorem *}, our third step is to deliver an upper bound of $P_{m}$ uniformly over $m$ and the main tool is also Lemma \ref{lemma 3}. There are only two different points. Firstly, the minimum distance between any two different $A_{jm}$ becomes $p_{N}$ instead of $p_{N}-2r$. Secondly, $N_{jm}^{r}$ is replaced by $N_{jm}$. Thus, similar to what we did above, for each $m$ and $N$(or relatively $Q$), we can construct a group of mutually independent random variables $\{W_{jm}:j=1,2dots, Q\}$ such that $W_{jm}$ is identical to $G_{jm}$. Meanwhile, we have, for any given $t>0$, 
\begin{align*}
    P_{m}\leq \mathbb{P}\left(\left|\sum_{j=1}^{Q}W_{jm}\right|>N\zeta-\sum_{j=1}^{Q}\xi_{jm}\right)+\sum_{j=2}^{Q}\mathbb{P}\left(\left|W_{jm}-G_{jm}\right|>\xi_{jm}\right).
\end{align*}
Note that, for any $t>0$, by letting $\xi_{jm}=\frac{1}{2}N_{jm}\zeta$, we have
\begin{align*}
    \mathbb{P}\left(\left|W_{jm}-G_{jm}\right|>\xi_{jm}\right)\leq 11\sqrt{2}(\frac{A}{\zeta}+1)^{1/2}(N_{jm}+N_{j-1,m})^{\tau}\psi_{\epsilon}(p_{N}),
\end{align*}
which yields
\begin{align*}
   \sum_{j=2}^{Q}\mathbb{P}\left(\left|W_{jm}-G_{jm}\right|>\xi_{jm}\right)\leq 2^{\tau+0.5}11(\frac{A}{\zeta}+1)^{0.5}M(Q,N,\tau)\psi_{\epsilon}(p_N).
\end{align*}
When there exists some universal $C>0$ such that $\sigma_{1}\leq \frac{Ct}{2^{d_2}}$, by setting $\xi_{jm}=\frac{1}{2C}N_{jm}\sigma_1$, we obtain
\begin{align*}
     \sum_{j=2}^{Q}\mathbb{P}\left(\left|W_{jm}-G_{jm}\right|>\xi_{jm}\right)\leq 22C^{1/3}M'(Q,N,\tau)\psi^{\frac{2}{3}}_{\epsilon}(p_N).
\end{align*}
Here the terms $M(Q,N,\tau)$ and $M'(Q,N,\tau)$ are directly inherited from the proof of Theorem \ref{Theorem *}. 
\par The final step is the application of the standard Bernstein inequality for independent data. This leads us to
\begin{align*}
    &\mathbb{P}\left(\left|\sum_{j=1}^{Q}W_{jm}\right|>N\zeta-\sum_{j=1}^{Q}\xi_{jm}\right)\leq \mathbb{P}\left(\left|\sum_{j=1}^{Q}W_{jm}\right|>\frac{1}{2}N\zeta\right)\\
    \leq& 2\exp\left(\frac{(N\zeta)^{2}}{8(\sum_{j=1}^{Q}Var(W_{jm})+\frac{1}{3}(\max_{j}N_{jm})AN\zeta)}\right)\\
    \leq&2\exp\left(\frac{Nd_{0}^{d}t^{2}}{2^{d_2+3}(\frac{Q}{N}d_{0}^{d}\Omega+\frac{1}{6}C_{2}p_{N}^{d_2}At)}\right),
\end{align*}
where $C_{2}$ is borrowed from the proof of Theorem \ref{Theorem *}. Then, we finish the proof.\\
\par \textbf{Proof of Corollary \ref{corollary 4}}
\par Our proof is motivated by the proof of Theorem A.2 in \citeA{xu2018sieve}. But we still demonstrate the whole proof here for self-completeness. 
\par We introduce the following truncation function. For some $A>0$, let
\begin{align*}
    f_{A}(x)=x1[|x|\leq A]+A1[x>A]-A1[x<-A].
\end{align*}
It is obvious that $||f_A||_{\infty}\leq A$ and $|f_{A}(Z_{s_{iN}})-E[f_{A}(Z_{s_{iN}})]|\leq 2A$. Furthermore, a key observation here is that $f_{A}$ is a Lipschitz continuous function on real line. More specifically, $|f_{A}(x)-f_A(x')|\leq |x-x'|$ holds for any $x\neq x'\in\mathbf{R}$. Then, according to Proposition 4, we know $\{f_{A}(Z_{s_{iN}}):s_{iN}\in\Gamma_N\}$ is also an array of random variables satisfying $L^p$-NED condition on $\epsilon_{s_{iN}}$'s and the coefficient of $L^p$-NED condition does not change up to a constant factor.
\par Then, we have
\begin{align*}
    \mathbb{P}(|N^{-1}S_{N}|>t)\leq& \mathbb{P}\Big(\Big|N^{-1}\sum_{i=1}^{N}(f_{A}(Z_{s_{iN}})-E[f_{A}(Z_{s_{iN}})])\Big|>t\Big)\\
    &+\mathbb{P}\Big(\Big|N^{-1}\sum_{i=1}^{N}(f^c_{A}(Z_{s_{iN}})-E[f^c_{A}(Z_{s_{iN}})])\Big|>t\Big),
\end{align*}
where $f^c_{A}(Z_{s_{iN}})=Z_{s_{iN}}-f_{A}(Z_{s_{iN}})$. Note that
\begin{align*}
    &\mathbb{P}\Big(\Big|N^{-1}\sum_{i=1}^{N}(f^c_{A}(Z_{s_{iN}})-E[f^c_{A}(Z_{s_{iN}})])\Big|>t\Big)\\
    \leq& \frac{2\max_{i,N}E|f^c_{A}(Z_{s_{iN}})|}{t}= \frac{2\max_{i,N}E|Z_{s_{iN}}1[|Z_{s_{iN}}|> A]|}{t}.
\end{align*}
According to Assumption \ref{assumption exponential}, we know $Z_{s_{iN}}$ has arbitrary order of moments. Hence, for any given $i$ and $N$ and $p,\ q\geq 1$ such that $\frac{1}{p}+\frac{1}{q}=1$, by applying Hölder's inequality, there exists some universal constant $C_q$ associated with $q$ such that
\begin{align*}
    E|f^c_{A}(Z_{s_{iN}})|\leq ||Z_{s_{iN}}||_{p}||1[|Z_{s_{iN}}|>A]||_{q}\leq C_q||Z_{s_{iN}}||_{p} \exp(-\frac{A}{2q\eta}).
\end{align*}
 Now, for any given universal $\Theta'>0$, by letting $A=\frac{\Theta'}{2q\eta}\log N$, when $\frac{\Theta'}{2q\eta}\log N\geq \frac{\kappa^2}{\eta}$,
 \begin{align*}
     \mathbb{P}\Big(\Big|N^{-1}\sum_{i=1}^{N}(f^c_{A}(Z_{s_{iN}})-E[f^c_{A}(Z_{s_{iN}})])\Big|>t\Big)\leq \frac{2\max_{i,N}||Z_{s_{iN}}||_{p}}{N^{\Theta'}t}.
 \end{align*}
 Meanwhile, together with Theorem \ref{Theorem *}, we have
 \begin{align*}
     &\mathbb{P}\Big(\Big|N^{-1}\sum_{i=1}^{N}(f_{A}(Z_{s_{iN}})-E[f_{A}(Z_{s_{iN}})])\Big|>t\Big)\\
     \leq& 2^{d_2+1}\exp\left(-\frac{Nd_{0}^{d}t^{2}}{C_{3}(\mathbf{V}d_{0}^{d}+C_{2}(p_N)^{d_2}\frac{\Theta'}{2q\eta}(\log N)t)}\right)\notag \\
    &+2^{d_2+\tau+\frac{1}{2}}11\left(\frac{2^{d_2}\Theta'\log N}{q\eta t}+1\right)^{1/2}M(Q,N,\tau)\psi_{\epsilon}(\frac{p_{N}}{3})\notag \\
    &+2^{d_2}\left(\frac{2^{d_2+2}A(N)\psi_{Z}(\frac{p_{N}}{3})}{t}\right)^{p}.
 \end{align*}
 
\section{Proof Related to Section 4}
\begin{lemma}(\textbf{Union Bound Argument})
\label{lemma 1 CV}
Define $\Theta$ as a compact subset of $\mathbf{R}^{m}$ with radius $R$. $\mathcal{F}$ and $\mathcal{G}$ are two sets of functions. Suppose the following three conditions about $\mathcal{F}$ and $\mathcal{G}$ are satisfied.
\begin{itemize}
    \item [1]$\mathcal{F}$ can be rewritten as $\{f_{\theta}:\theta\in\Theta\}$.
    \item [2]There exists an injective mapping from $\Theta$ to $\mathcal{F}$, denoted as $T:\theta\rightarrow f_{\theta}$.
    \item [3]Given a fixed $\theta$ and its neighbor $B(\theta,\epsilon_{n})$, there exists an $n_1\in\mathbf{N}^{+}$ such that for any $n>n_1$, when $\theta\in\Theta$, $\theta\neq \theta '$ and $||\theta-\theta'||_{E}\leq \epsilon_{n}$, there exists a function in $\mathcal{G}$, denoted as $G_{\theta '}$, and a constant $M>0$ satisfying 
\begin{align*}
    ||f_{\theta}-f_{\theta '}||_{\infty}\leq M ||\theta-\theta'||_{E} G_{\theta '},
\end{align*}
where $||\cdot||_{E}$ denotes the Euclidean norm.
\end{itemize}
We further define $\{U_{in}\}$ is an array of random variables such $$\sup_{G\in\mathcal{G}}|\sum_{i=1}^{N}E[G(X_i)U_{in}]|:=M_{G}<+\infty$$ holds for every $n$. Then, for any $t>0$ satisfying $\frac{t}{6 MM_G}\leq \epsilon_{n}$, we have
\begin{align*}
    &\mathbb{P}(\sup_{f\in\mathcal{F}}|P_{n}(fU_{n})-P(fU_{n})|>t)\\
    \lesssim &(\frac{6\rho MM_{G}}{t})^{m}\max_{f}\mathbb{P}(|P_{n}(fU_{n})-P(fU_{n})|>\frac{t}{3})\\
    +&(\frac{6\rho MM_{G}}{t})^{m}\max_{G}\mathbb{P}(|P_{n}(GU_{n})-P(GU_{n})|>2M_{G}),
\end{align*}
where $P_{n}(gU_n)=\frac{1}{n}\sum_{i=1}^{n}g(X_{i}U_{in})$ and $P(gU_n)=\frac{1}{n}\sum_{i=1}^{n}E[g(X_{i})U_{in}]$, for a given function $g$ and sample $\{X_{i},i=1,\dots,n\}$.
\end{lemma}
\begin{lemma}
    \label{lemma 4-1}
    Assume that Assumptions \ref{Asp 22}-\ref{ASP 10}, by letting $h=(\log N/N)^{\frac{1}{2\alpha+D}}$, where $\alpha$ is from Assumption \ref{ASP 9}, we have 
 \begin{align*}
     \sup_{x\in\mathcal{B}_{D}(\rho_N)}\left|\frac{1}{Nh^{D}}\sum_{i=1}^{N}\left(K\left(\frac{X_{s_{iN}-x}}{h}\right)Y_{s_{iN}}-E\left[K\left(\frac{X_{s_{iN}-x}}{h}\right)Y_{s_{iN}}\right]\right)\right|\lesssim \mathcal{O}_{\mathbb{P}}\left(\frac{\log N}{N}\right)^{\frac{\alpha}{2\alpha+D}},
 \end{align*}    
 where $\mathcal{B}_{D}(R_N)$ is defined in Section 4. 
\end{lemma}
\par \textbf{Proof of Lemma \ref{lemma 1 CV}}
\par The proof of lemma \ref{lemma 1 CV} is very simple and straight. Firstly, please note that, due to the compactness of $\Theta$. For any 
$\xi>0$, its $\xi$-covering number with respect to Euclidean norm, denoted as $\mathcal{N}(\xi,\Theta,||\cdot||_{E})$, can be upper bounded as follow,
\begin{align*}
    \mathcal{N}(\xi,\Theta,||\cdot||_{E})\lesssim (\frac{\rho}{\xi})^{m}.
\end{align*}
Suppose $\Theta_{S}:=\{\theta_{1},\dots,\theta_{S}\}$ is the smallest $\frac{t}{6MM_{G}}$-net of set $\Theta$ with respect to Euclidean norm $||\cdot||_{E}$. According to conditions 1 and 2, we have 
\begin{align*}
    &\sup_{f}|P_{n}(fU_{n})-P(fU_{n})|=\sup_{\theta}|P_{n}(f_{\theta}U_{n})-P(f_{\theta}U_{n}))|\\
    \leq &\sup_{\theta,\theta^{k}} |P_{n}((f_{\theta}-f_{\theta^{k}})U_{n})|+\sup_{\theta,\theta^{k}}|P((f_{\theta}-f_{\theta^{k}})U_{n})|+\max_{\theta^{k}\in\Theta_{S}}|P_{n}(f_{\theta^{k}}U_{n})-P(f_{\theta^{k}}U_{n})|.
\end{align*}
According to condition 3 and the assumption that $\frac{t}{6MM_G}\leq \epsilon_{n}$, when $n\geq n_1$, there exists a $G_{\theta^{k}}\in\mathcal{G}$, such that
\begin{align*}
    ||f_{\theta}-f_{\theta^{k}}||_{\infty}\leq M ||\theta-\theta'||_{E} G_{\theta^{k}}, 
\end{align*}
which indicates that, for any $n\geq n_0\lor n_1$,
\begin{align*}
    &\sup_{f}|P_{n}(fU_{n})-P(fU_{n})|\\
    \leq& \frac{t}{6M_{G}}\max_{\theta_{k}}|P_{n}(G_{\theta_{k}}U_{n})-P(G_{\theta_{k}}U_{n})|)+\frac{t}{3}+\max_{\theta^{k}\in\Theta_{S}}|P_{n}(f_{\theta^{k}}U_{n})-P(f_{\theta^{k}}U_{n})|.
\end{align*}
Therefore, we have
\begin{align*}
   &\mathbb{P}(\sup_{f}|P_{n}(fU_{n})-P(fU_{n})|>t)\\
   \leq &\mathbb{P}(\max_{\theta^{k}}\frac{t}{6M_{G}}|P_{n}(G_{\theta^{k}}U_{n})-P(G_{\theta^{k}}U_{n})|+\frac{t}{3}+\max_{\theta^{k}}|P_{n}(f_{\theta^{k}}U_{n})-P(f_{\theta^{k}}U_{n})|>t)\\
   \leq &S\max_{\theta^{k}}\mathbb{P}(\frac{t}{6M_{G}}|P_{n}(G_{\theta^{k}}U_{n})-P(G_{\theta^{k}}U_{n})|>\frac{t}{3})+S\max_{\theta^{k}}\mathbb{P}(|P_{n}(f_{\theta^{k}}U_{n})-P(f_{\theta^{k}})|>\frac{t}{3})\\
   \lesssim &(\frac{6\rho MM_{G}}{t})^{m}\max_{f}\mathbb{P}(|P_{n}(fU_{n})-P(fU_{n})|>\frac{t}{3})\\
   &+(\frac{6\rho MM_{G}}{t})^{m}\max_{G}\mathbb{P}(|P_{n}(GU_{n})-P(GU_{n})|>2M_{G}).
\end{align*}
\textbf{Q.E.D}\\
\par \textbf{Proof of Lemma \ref{lemma 4-1}}
\par \textbf{Step 1}(Uniform Bound Argument)
\par To simplify our notation, here we use $i$ to represent subscript $s_{iN}$. That is to say, we denote $(X_{s_{iN}},Y_{s_{iN}})$ as $(X_i,Y_i)$ for simplicity. We first truncate the output $Y_{i}$'s based on the following decomposition. Given a number sequence $\{A_{N}\}$ such that $\lim_{N}A_{N}\nearrow +\infty$, we have $Y_i=Y_{i1}+Y_{i2}$, where $Y_{i1}=Y_{i}1[|Y_{i}|\leq A_{N}]$ and $Y_{i2}=Y_{i}1[|Y_i|>A_{N}]$. Then,
\begin{align*}
    &\sup_{x\in\mathcal{X}}\left|\frac{1}{Nh^{D}}\sum_{i=1}^{N}\left(K(\frac{X_{i}-x}{h})Y_i-E[K(\frac{X_{i}-x}{h})Y_i]\right)\right|\\
    \leq &\sup_{x\in\mathcal{X}}\left|\frac{1}{Nh^{D}}\sum_{i=1}^{N}K(\frac{X_{i}-x}{h})Y_{i1}-E[K(\frac{X_{i}-x}{h})Y_{i1}]\right|\\
    + &\sup_{x\in\mathcal{X}}\left|\frac{1}{Nh^{D}}\sum_{i=1}^{N}K(\frac{X_{i}-x}{h})Y_{i2}\right|+\sup_{x\in\mathcal{X}}\left|\frac{1}{Nh^{D}}\sum_{i=1}^{N}E[K(\frac{X_{i}-x}{h})Y_{i2}]\right|\\
    := \mathbf{S}_1+\mathbf{S}_2+\mathbf{S}_3.
\end{align*}
\par Now we focus on bounding $\mathbf{S}_2$ and $\mathbf{S}_3$.
\par It is obvious that, for any given $C>0$ and $a_{N}>0$, we have 
\begin{align*}
    \mathbb{P}(\mathbf{S}_{3}>Ca_{N})\leq \mathbb{P}(\bigcup_{i=1}^{N}\{|Y_i|>A_{N}\})\leq \sum_{i=1}^{N}\mathbb{P}(|Y_i|>A_{N})\leq \frac{N\max_{i}E|Y_i|^{p_0}}{A_{N}^{p_0}}.
\end{align*}
Meanwhile, based on some simple algebra, we can get, for any $i$,
\begin{align*}
    &h^{-D}\left|E[K(\frac{X_{i}-x}{h})Y_{i2}]]\right|\leq A_{N}^{-(p_{0}-1)}h^{-D}E\left|K(\frac{X_{i}-x}{h})Y_{i2}^{p_0}\right|\\
    =&A_{N}^{-(p_0-1)}h^{-D}\int_{\mathbf{R}^{D}}K(\frac{z-x}{h})E[|Y_i|^{p_0}|X_i=z]f_{i}(z)dz\\
   =& A_{N}^{-(p_0-1)}\int_{\mathbf{R}^{D}}K(u)E[|Y_i|^{p_0}|X_i=x+uh]f_{i}(x+uh)du\\
   \leq &A_{N}^{-(p_0-1)}(\sup_{x}E[|Y_i|^{p_0}|X_i=x]f_{i}(x))\int_{\mathbf{R}^{D}}K(u)du\lesssim A_{N}^{-(p_0-1)}.
\end{align*}
This leads to $\mathbf{S}_2\lesssim A_{N}^{-(p_0-1)}$. Thus, given the $a_{N}$ mentioned above, if we want both $\mathbf{S}_2$ and $\mathbf{S}_3$ convergence to $0$ at some rates no slower than $a_{N}$, we have to ensure
\begin{align}
    \label{P4.2-1}
    N/A_{N}^{p_0}=o(1)\ \text{and}\ A_{N}^{-(p_0-1)}\lesssim a_{N}.
\end{align}
\par Now we focus on dealing with $\mathbf{S}_1$. Denote $f(z,x,h)=K(\frac{X_{i}-x}{h})$ and $\mathcal{F}=\{f(z,x,h):x\in\mathcal{X}\}$. Then, for some user defined constant $C_{a}>0$,
\begin{align*}
    \mathbb{P}(\mathbf{S}_1>C_{a}a_{N})=\mathbb{P}\left(\sup_{f\in\mathcal{F}}\left|\frac{1}{N}\sum_{i=1}^{N}\left(f(X_i,x,h)Y_{i1}-E[f(X_i,x,h)Y_{i1}]\right)\right|>C_{a}a_{N}h^{D}\right).
\end{align*}
Based on the Lipschitz property of kernel function, given any $x\in\mathcal{X}$ and $x'\in B(x,h)$, there exists some universal constant $c_D>0$ such that
\begin{align}
    \label{P4.2-2}
    ||f(z,x,h)-f(z,x',h)||_{\infty}&\leq c_Dh^{-1}||x-x'||_{E}\prod_{k=1}^{D}1[|z^{k}-x'^{k}|\leq 2h]\notag \\
    &:=M||x-x'||_{E}G(z,x',h).
\end{align}
Thus, we construct set $\mathcal{G}=\{G(z,x,h):x\in\mathcal{X}\}$. By denote the Lebesgue measure on $\mathbf{R}^{D}$ as $Leb$, We can easily obtain that 
\begin{align*}
    \sup_{G\in\mathcal{G}}|\frac{1}{N}\sum_{i=1}^{N}E[G(Z_{i},x,h)Y_{i1}]|\leq Leb(B(x,h))\max_{i}|E[Y_{i}]|\leq Ch^{D}.
\end{align*}
Hence, based on \eqref{P4.2-1} and \eqref{P4.2-2}, the $M$ and $M_{G}$ in Lemma \ref{lemma 1 CV} are equal to $c_{D}h^{-1}$ and $Ch^{D}$ respectively. To apply Lemma \ref{lemma 1 CV}, we only need to check whether $\frac{a_{N}h^{D}}{6c_{D}h^{-1}Ch^{D}}\leq h$ holds for sufficiently large $n$. Considering that $a_{N}=o(1)$, this is true. Thus, according to Lemma \ref{lemma 1 CV}, we obtain that
\begin{align*}
    \label{P4.2-3}
     &\mathbb{P}(\mathbf{S}_1>C_{a}a_{N})\notag\\
     =&\mathbb{P}\left(\sup_{f\in\mathcal{F}}\left|\frac{1}{N}\sum_{i=1}^{N}\left(f(X_i,x,h)Y_{i1}-E[f(X_i,x,h)Y_{i1}]\right)\right|>C_{a}a_{N}h^{D}\right)\notag\\
     \lesssim &(\frac{\rho_N}{a_{N}h})^{D}\sup_{x}\mathbb{P}\left(\left|\frac{1}{N}\sum_{i=1}^{N}\left(f(X_i,x,h)Y_{i1}-E[f(X_i,x,h)Y_{i1}]\right)\right|>C_{a}a_{N}h^{D}\right)\\
     &+(\frac{\rho_N}{a_{N}h})^{D}\sup_{x}\mathbb{P}\left(\left|\frac{1}{N}\sum_{i=1}^{N}\left(G(X_i,x,h)Y_{i1}-E[G(X_i,x,h)Y_{i1}]\right)\right|>2Ch^{D}\right)\\
     :=& Q_1+Q_2.
\end{align*}
Obviously, it can be easily figured out that $Q_{2}$ is a higher order term of $Q_1$. Hence, we now only focus on $Q_1$.\\
\par \textbf{Step 2}(Application of Corollary \ref{Theorem 2})
\par Let $W_{s_{iN}}=f(X_{s_{iN}},x,h)Y_{s_{iN}1}=K(\frac{X_{s_{iN}}-x}{h})Y_{{s_{iN}}1}$, where $X_{s_{iN}1}$ and $Y_{s_{iN}1}$ are equal to $X_{i1}$ and $Y_{i1}$ defined above. We instantly know that, together with the assumptions about NED conditions mentioned Section 4.2, we have, for $p\geq 1$,
\begin{align*}
    ||W_{s_{iN}}-E[W_{s_{iN}}|\mathcal{F}_{iN}(r)]||_{p}\leq 2DLip(K)h^{-1}r^{-\gamma}\lesssim h^{-1}r^{-\gamma}.
\end{align*}
Thus, the $A(N)$ for array $\{W_{s_{iN}}\}$ is of order $h^{-1}$. Similarly, we still denote $W_{s_{iN}}$ as $W_{i}$ for simplicity. 
\par Now, we need to specify $\mathbf{V}_{2}$ and $\max_{i,N}||W_i||_{\infty}$ so that Corollary \ref{Theorem 2} could be applied. Obviously $\max_{i,N}||W_i||_{\infty}\leq ||K||_{\infty}^{D}A_{N}$. Now we start to calculate the upper bound of $\mathbf{V2}$. First, please note that, uniformly over $i$, according to the M2 in Assumption 8, we have 
\begin{align*}
    Var(W_{i})\lesssim h^{D}\ \text{and}\ ||W_{i}||_{p_0}^{2}\lesssim h^{\frac{2D}{p_0}}.
\end{align*}
Then, note that $p_0$ here is actually $2+\delta$ in Proposition \ref{prop *}. Then, together with the assumption that $d_0=1$, by setting $\tau_{N}=h^{-\frac{D}{d_2}}$, we have
\begin{align*}
  \mathbf{V2}\lesssim h^{D}+h^{-D}\Sigma+h^{\frac{2\gamma D}{d_2}-D-2}+h^{\frac{D\gamma}{d_2}-0.5D-1}+h^{\frac{D}{d_2}(\frac{(p_0-2)\gamma}{p_0}-d_2)+\frac{2D}{p_0}}.
\end{align*}
Furthermore, for any given $x$ and $i\neq j$, we have
\begin{align*}
    &E|W_{i}W_{j}|\\
    \leq& \int_{\mathbf{R}^{2D}}K(\frac{X_i-x}{h})K(\frac{X_j-x}{h})E[|Y_iY_j||X_i=u,X_j=v]f_{(X_i,X_j)}(u,v)dudv\\
    \leq& h^{2D}\int_{\mathbf{R}^{2D}}K(u)K(v)E[|Y_iY_j||X_i=x+uh,X_j=x+vh]f_{(X_i,X_j)}(x+uh,x+vh)dudv\\
    \leq& \left(\max_{i\neq j}\sup_{s,t} E[|Y_iY_j||X_i=s,X_j=t]f_{(X_i,X_j)}(s,t)\right)h^{2D}||K||_{1}^{2}.
\end{align*}
Thus, based on the definition of $\Sigma$, we can show $\Sigma\lesssim h^{2D}$. According to R4 in Assumption \ref{ASP 10}, we immediately know there exists some universal constant $C_{V2}>0$ such that $\mathbf{V2}\leq C_{V2}h^{D}$. 
\par Now we define $h=(\frac{\log N}{N})^{\frac{1}{2\alpha+D}}$, $a_{N}=\sqrt{\frac{\log N}{Nh^{D}}}=(\frac{\log N}{N})^{\frac{\alpha}{2\alpha+D}}$, $A_{N}=N^{\frac{1}{p_0}}\xi_{N}$ and $p_N=(a_N A_N)^{-\frac{1}{d_2}}$, where $\{\xi_{N}\}$ is arbitrary sequence diverging to $+\infty$ at a rate no faster than $\log N$. For some user defined $C_a>0$, according to Corollary \ref{Theorem 2}, we have
\begin{align*}
    Q_1&\lesssim (\frac{\rho_N}{a_{N}h})^{D}\exp\left(-\frac{C_{a}^{2}Na^{2}_{N}h^{2D}}{C_{V2}h^{D}+C_2||K||_{\infty}^{D}(a_N A_N)^{-1}A_Na_N h^{D}}\right)\\
    &+ (\frac{\rho_N}{a_{N}h})^{D}(\frac{A_N}{a_{N}h^{D}})^{0.5}(a_{N}A_{N})^{\gamma/d_2} +(\frac{\rho_N}{a_{N}h})^{D}(\frac{A_N}{a_{N}h^{D}})^{p}(a_{N}A_{N})^{\frac{\gamma}{d_2} p}\\
    &:= Q_{11}+Q_{12}+Q_{13}.
\end{align*}
Recall $a_{N}=(\log N/N)^{\frac{\alpha}{2\alpha+D}}$, $h=(\log N/N)^{\frac{1}{2\alpha+D}}$ and $A_{N}=N^{1/p_0}\xi_{N}$. It is obvious that, by choosing sufficiently large $C_{a}$, we can have $Q_{11}=o(1)$. Now we only need to find sufficient conditions such that $Q_{12}\lor Q_{13}=o(1)$. Actually, based on the definitions of $a_N$, $h$ and $A_N$ and the fact that $\xi_{N}/\log N=O(1)$, we only need to consider two terms $N$, $\log N$ and $\xi_{N}$. Thus, it suffices to focus only on the exponent of $N$. Then, according to some direct algebra, it can be shown that the sufficient condition for $Q_{12}=o(1)$ is the R2 of Assumption \ref{ASP 10} and the sufficient condition for $Q_{13}=o(1)$ is R3 of Assumption \ref{ASP 10}. Then we finish the proof. \textbf{Q.E.D}\\
\par \textbf{Proof of Theorem \ref{Theorem 4.2}}
\par Within this proof, we always admit that $h=(\frac{\log N}{N})^{\frac{1}{2\alpha+D}}$. Proof of Theorem \ref{Theorem 4.2} is almost the same as the proof of Theorem 3 in \citeA{jenish2012nonparametric}. Hence we only highlight the difference here. A basic decomposition of $\hat{\theta}(x)-\theta(x)$ is as follow,
\begin{align*}
    \hat{\theta}(x)-\theta(x)=\begin{pmatrix}
    \hat{m}(x)-m(x)\\ (\triangledown \hat{m}(x)-\triangledown m(x))h
    \end{pmatrix}
    =S_{N}^{-1}(x)\underbrace{\left\{T_N(x)-S_{N}(x)\begin{pmatrix}
    m(x)\\ \triangledown m(x)h
    \end{pmatrix}\right\}}_{:=W_N (x)},
\end{align*}
where, for $k=1,\dots,D+1$, by denoting $v^{k}$ as the k-th element of vector $v$ and location $s_{iN}$ as $iN$,
\begin{align*}
    &W_{N}^{k}(x)=\frac{1}{Nh^{D}}\sum_{i=1}^{N}K\left(\frac{X_{iN}-x}{h}\right)U_{i}^{k}(x)\omega_{iN}\\
    &\omega_{iN}= Y_{iN}-m(x)-\triangledown m(x)^{T}(X_{iN}-x).
\end{align*}
$\triangledown m(x)$ denotes the gratitude of $m(x)$. Hence, it suffices to finish the proof if we demonstrate the following three things,
\begin{enumerate}[(1)]
    \item $\sup_{x\in\mathcal{X}}||U^{-1}(x)||=O_{p}(1)$,
    \item $\sup_{x\in\mathcal{X}}|E[W_{N}^{1}(x)]|=O_{p}(h^{\alpha})$,
    \item $\sup_{x\in\mathcal{X}}|W_{N}^{1}(x)-E[W_{N}^{1}(x)]|=\sqrt{\frac{\log N}{Nh^{D}}}$,
\end{enumerate}
where $||A||$ means Euclidean norm of matrix A. Furthermore, we have the following decomposition
\begin{align*}
    &\sup_{x\in\mathcal{X}}|W_{N}^{1}(x)-E[W_{N}^{1}(x)]|\\
    \leq & \sup_{x\in\mathcal{X}}\left|\frac{1}{Nh^{D}}\sum_{i=1}^{N}\left(K\left(\frac{X_{iN}-x}{h}\right)U_{i}^{1}(x)Y_{iN}-E\left[K\left(\frac{X_{iN}-x}{h}\right)U_{i}^{1}(x)Y_{iN}\right]\right)\right|\\ + &\sup_{x\in\mathcal{X}}||m||_{\infty}\left|\frac{1}{Nh^{D}}\sum_{i=1}^{N}\left(K\left(\frac{X_{iN}-x}{h}\right)U_{i}^{1}(x)-E\left[K\left(\frac{X_{iN}-x}{h}\right)U_{i}^{1}(x)\right]\right)\right|\\ + &\sup_{x\in\mathcal{X}}||\triangledown m||\left\|\frac{1}{Nh^{D}}\sum_{i=1}^{N}\left((X_{iN}-x)K\left(\frac{X_{iN}-x}{h}\right)U_{i}^{1}(x)-E\left[(X_{iN}-x)K\left(\frac{X_{iN}-x}{h}\right)U_{i}^{1}(x)\right]\right)\right\|\\
    := & S_1 +S_2 +S_3,
\end{align*}
where $U_{i}^{1}(x)=1$. According to Assumption \ref{ASP 9}, we instantly know that $||m||_{\infty}\lor ||\triangledown m||<+\infty$ holds universally. Therefore, due to the compactness of support set $\mathcal{X}$, $S_2$ and $S_3$ enjoys uniform boundedness. Hence, it suffices to only consider $S_1$ since $S_2$ and $S_3$ are special cases of $S_1$. According to Lemma \ref{lemma 4-1}, we can show that, under Assumptions \ref{ASP 7}, \ref{ASP 8} and \ref{ASP 10}, $S_{1}\lesssim (\log N/N)^{\frac{\alpha}{2\alpha+D}}$. 
\par Meanwhile, together with Assumption \ref{Asp 22} of kernel function and Assumption \ref{ASP 9}, (2) can also be obtained based on standard Taylor expansion argument. Thus it is omitted here. To prove (1), it is equivalent to prove matrix $U(x)^{-1}$ is asymptotically existed, which can be done by using Proposition \ref{prop *}. Since procedure of this is almost the same as the proof of Lemma 1 in \citeauthor{jenish2012nonparametric}, it is omitted as well. Then, we finish the proof.

\end{document}